\documentclass[12pt,twoside]{article}

\usepackage{amsmath}
\usepackage{amssymb}
\usepackage{amsthm}
\usepackage{amscd}
\usepackage{color}
\usepackage{accents}
\usepackage{mathrsfs}
\usepackage{url}

\usepackage{enumitem}
\usepackage{euscript}
\usepackage{multirow}
\usepackage{latexsym,amsbsy}
\usepackage[pdftex]{graphicx}
\usepackage{epsfig}
\usepackage{subdepth}

\usepackage{footnote}

\usepackage{float}
\usepackage{color}
\usepackage{relsize}
\usepackage[geometry]{ifsym}
\usepackage{pdfpages}
\usepackage[makeroom]{cancel}
\usepackage[normalem]{ulem}

\usepackage{eqnarray}

\usepackage[linktoc=all,bookmarks,hypertexnames=false,debug,pagebackref=false,colorlinks=true,
linkcolor=blue,citecolor=blue,urlcolor=blue]{hyperref}
\usepackage{bookmark}

\usepackage{natbib}

\usepackage{apalike}
\makeatletter
\def\@cite#1#2{#1\if@tempswa , #2\fi}
\makeatother

\usepackage{tocloft}
\usepackage[multiple]{footmisc}

\numberwithin{equation}{section}

\theoremstyle{plain}
\newtheorem{theorem}{Theorem}[section]

\newtheorem{lemma}[theorem]{Lemma}
\newtheorem{corollary}[theorem]{Corollary}
\newtheorem{remark}[theorem]{Remark}
\newtheorem{example}[theorem]{Example}

\makeatletter
\renewcommand*{\@fnsymbol}[1]{\@arabic{#1}}
\newcommand*\bigcdot{\mathpalette\bigcdot@{.8}}
\newcommand*\bigcdot@[2]{\mathbin{\vcenter{\hbox{\scalebox{#2}{{\hskip 1pt}$\m@th#1\bullet$}}}}}
\makeatother

\def\E {\mathbb{E}\hskip1pt}
\def\N {\mathbb{N}}
\def\P {\mathbb{P}}
\def\R {\mathbb{R}}
\def\Var {\mathrm{Var}\hskip1pt}
\def\Cov {\mathrm{Cov}}

\newcommand{\dd}{{{\hskip2pt} \mathrm{d}}}
\def\i{\mathrm{i}\hskip1pt}

\newcommand{\mcalB}{\mathcal{B}}
\newcommand{\mcalC}{\mathcal{C}}
\newcommand{\mcalF}{\mathcal{F}}
\newcommand{\mcalN}{\mathcal{N}}
\newcommand{\mcalR}{\mathcal{R}}
\newcommand{\mcalS}{\mathcal{S}}

\newcommand{\mcalW}{\mathcal{W}}
\newcommand{\mcalX}{\mathcal{X}}
\newcommand{\mcalY}{\mathcal{Y}}
\newcommand{\mcalZ}{\mathcal{Z}}

\newcommand{\wtilX}{\widetilde{X}}

\def\cid{\stackrel{\mathcal{L}}\longrightarrow}
\def\cip{\stackrel{\mathcal{P}}\longrightarrow}
\def\cwip{\stackrel{w\mathcal{P}}\longrightarrow}

\def\diag{\mathop{\mathrm{diag}}}
\def\sinc{\mathop{\mathrm{sinc}}}
\def\tr{\mathop{\mathrm{tr}}}

\def\goesto#1#2{{\hskip 1pt} #1 {\hskip 1pt} \mapsto {\hskip 1pt} #2 {\hskip 1pt}}

\newcommand{\wtilde}{\widetilde}
\newcommand{\eqdist}{\overset{\mathcal{L}}{=}}

\newcommand{\bfone}{\boldsymbol{1}}

\begin{document}

\title{\vspace{-35pt}
\Large\textbf{Random Linear Modulation with Spherically Symmetric Modulators}}

\author{{Armine Bagyan}$^{1,2}$ {and Donald Richards}$^{1,3}$}

\date{ }

\maketitle

\phantom{a}
\vspace{-25pt}

\begin{abstract}
We consider the modulation of data given by random vectors $X_n \in \R^{d_n}$, $n \in \N$.  For each $X_n$, one chooses an independent modulating random vector $\Xi_n \in \R^{d_n}$ and forms the projection $Y_n = \Xi_n'X_n$.  It is shown, under regularity conditions on $X_n$ and $\Xi_n$, that $Y_n|\Xi_n$ converges weakly in probability to a normal distribution.  More broadly, the conditional joint distribution of a family of projections constructed from random samples from $X_n$ and $\Xi_n$ is shown to converge weakly to a matrix normal distribution.  We derive, ${via}$ G.~P\'olya's characterization of the normal distribution, a necessary and sufficient condition on $Y_n$ for $\Xi_n$ to be normally distributed, and we show that our results motivate generalizations of P\'olya's theorem.  When $\Xi_n$ has a spherically symmetric distribution we deduce, through I.~J.~Schoenberg's characterization of the spherically symmetric characteristic functions on Hilbert spaces, that the probability density function of $Y_n|\Xi_n$ converges pointwise in certain $p$th means to a mixture of normal densities and a rate of convergence is quantified, resulting in uniform convergence.  The cumulative distribution function of $Y_n|\Xi_n$ is shown to converge uniformly in those $p$th means to the distribution function of the same mixture, and a Lipschitz property is obtained.  Examples of distributions for $X_n$ that satisfy our results include the Bingham distributions on hyperspheres of random radii, uniform distributions on hyperspheres and hypercubes of random volumes, and multivariate normal distributions; and examples of such $\Xi_n$ include the multivariate $t$-, multivariate Laplace, and spherically symmetric stable distributions.  

\medskip



\smallskip
\noindent
{\em Running head}: Random Modulation with Spherical Symmetry.
\end{abstract}

\footnoterule
\footnotesize{
\noindent
$^1$Department of Statistics, Pennsylvania State University, University Park, PA 16802, U.S.A. \\
$^2$E-mail address: aub171@psu.edu \\
$^3$E-mail address: dsr11@psu.edu
}

\tableofcontents

\normalsize

\section{Introduction and motivation}
\label{sec_introduction}

Random modulation, in which several random signals are combined to form a new signal \citep{Black,Papoulis}, is well known from its role in amplitude modulation (AM) and frequency modulation (FM) radio broadcasting.  Random modulation is now applied widely, in fields such as electric power devices \citep{Blaabjerg_etal}, underwater ranging and detection \citep{Cochenour_etal}, autonomous vehicles \citep{Hwang}, radio-frequency identification (RFID) security \citep{Roy_etal}, atmospheric research \citep{She_etal}, medical technologies \citep{Tang_Clement}, wireless communications \citep{vanTrees}, and pathogen detection \citep{Yang_etal}.  

We are motivated here by questions arising from linear random modulation of high-dimensional data.  For each dimension $d_n$, $n=1,2,3,\ldots$, we are given datum in the form of a random vector $X_n \in \R^{d_n}$.  On choosing an independent modulating random vector $\Xi_n \in \R^{d_n}$, and forming the projection $Y_n = \Xi_n'X_n$, we study the limiting conditional distribution of $Y_n|\Xi_n$ under regularity assumptions on $X_n$ and $\Xi_n$, as $d_n \to \infty$.  

Linear modulation appears in mathematical statistics prominently in the study of low-dimensional projections of high-dimensional vectors, where the notable results of \cite{Eaton}, \cite{Diaconis_Freedman}, \cite{Huber}, and \cite{Duembgen_etal} have spawned an extensive literature.  Also noteworthy are \cite{Lok_Lehnert}, who studied linearly modulated communication systems; \cite{Loperfido}, in the area of detecting financial outliers; and \cite{Davidov_Peddada}, who formulated the theoretical foundations of ordered projections of multivariate data and gave applications to the analysis of toxicological and carcinogenic data.  

Among the cited literature, we emphasize the work of \cite{Duembgen_etal} and \cite{Wee_Tatikonda} who derived, along with other results, the weak convergence properties of the conditional distribution functions of $Y_n$.  Our results are also concerned with those conditional distributions.  However, we proceed using different methods that yield the convergence properties of both the conditional probability density and the conditional cumulative distribution functions of such projections.  

Our work is motivated proximately by \cite{Bagyan}, who derived $L^2$-pointwise convergence results for the conditional density and distribution functions of $Y_n|\Xi_n$ when $\Xi_n$ is normally distributed.  For absolutely continuous distributions, uniform convergence of the distribution functions follows from their weak convergence \citep[p.~62]{Zolotarev}; however it is generally more difficult to derive the convergence properties of the corresponding density functions.  Thus we extend the results of \cite{Bagyan} by obtaining $L^p$-pointwise convergence results for the density function of $Y_n|\Xi_n$ and $L^p$-uniform convergence results for its distribution function; also, we extend these results to the case in which $\Xi_n$ is spherically distributed.  

In Section \ref{sec_gaussian_case}, we suppose that $\Xi_n \sim \mcalN_{d_n}(0,I_{d_n})$, the multivariate standard normal distribution.  \citet[Section 2.6]{Bagyan} also studied this case with random sampling conducted on $X_n$ and $\Xi_n$, yielding data $X_{n,1},\ldots,X_{n,k}$, and $\Xi_{n,1},\ldots,\Xi_{n,l}$, respectively, and derived the limiting unconditional distribution of the collection of projections $\{\Xi_{n,j}'X_{n,r}, j=1,\ldots,l, r=1,\ldots,k\}$.  Throughout the article, all results are derived under the assumption that $X_n$ satisfies the regularity conditions \ref{condition_1} and \ref{condition_2}.  By adapting the approach of \cite{Bagyan}, we obtain in Theorem \ref{thm_resampling_Xi_n_X_n} the limiting weak distribution of this collection of projections, conditional on $\Xi_{n,1},\ldots,\Xi_{n,l}$.  Theorem \ref{thm_resampling_Xi_n_X_n} is related to numerous articles (cf., \cite{Diaconis_Freedman}, \cite{Duembgen_etal}, \cite{Bickel_etal}) that explain why many unit-length projections of a high-dimensional random vector are approximately normally distributed, and our proof motivates the results in Sections \ref{sec_Lp_conv_cond_pdf} and \ref{sec_Lp_conv_cond_cdf}.  

For the case in which $\Xi_n$ is spherically symmetric, we obtain in Theorem \ref{thm_modulation_charac_2} a necessary and sufficient condition for $\Xi_n$ to be normally distributed; this result, which may be a new uniqueness property of the multivariate normal distribution, will be derived \textit{via} a celebrated result of \cite{Polya} that characterizes the normal distribution through the distribution of linear functions of independent random variables.  We also show in Remark \ref{rem_C_1_prime_sec_2} that our results motivate generalizations of P\'olya's theorem.  

The data $X_n$ are assumed throughout to satisfy the regularity conditions \ref{condition_1} and \ref{condition_2}, so we provide in Section \ref{sec_examples} examples of distributions that satisfy those assumptions (cf.~\citet[Section 2]{Duembgen_etal} for other examples).  Our examples include dilated Bingham distributions on hyperspheres, uniform distributions on Euclidean balls and on hypercubes, and multivariate normal distributions.  Further it is shown that the multivariate $t$-distributions satisfy \ref{condition_2} but do not satisfy \ref{condition_1}.  

From Section \ref{sec_Lp_conv_cond_pdf} onwards, we assume that the modulating vector $\Xi_n$ is spherically symmetric and we consider the convergence properties of $f_{Y_n|\Xi_n}$, the probability density function of $Y_n|\Xi_n$.  By applying a famous theorem of \cite{Schoenberg}, which characterizes the class of spherically symmetric characteristic functions on Hilbert spaces, we derive conditions such that certain powers, $[f_{Y_n|\Xi_n}(\cdot)]^k$, $k \in \N$, converge $L^p$-pointwise-in-mean to corresponding powers, $[f_{\mcalN_1(0,\sigma^2 V^2)}(\cdot)]^k$, of a normal mixture density, where the random variable $V$ determined by $\Xi_n$.  Further we deduce the pointwise convergence in certain $p$th means of $f_{Y_n|\Xi_n}(\cdot)$ to $f_{\mcalN_1(0,\sigma^2 V^2)}(\cdot)$.  Motivated by results of \cite{Meckes}, and \cite{Wee_Tatikonda}, who obtained quantitative asymptotics for convergence results in projection analysis, we obtain an inequality for the difference, $|\E[f_{Y_n|\Xi_n}(\cdot)]^k - \E[f_{\mcalN_1(0,\sigma^2 V^2)}(\cdot)]^k|$, hence yielding a rate of convergence in terms of the regularity conditions \ref{condition_1} and \ref{condition_2}.  

In Section \ref{sec_Lp_conv_cond_cdf}, we provide conditions under which $k$th powers of $F_{Y_n|\Xi_n}(\cdot)$, the cumulative distribution function $Y_n|\Xi_n$, converge uniformly in mean to $k$th powers of the corresponding mixture distribution function $F_{\mcalN_1(0,\sigma^2 V^2)}(\cdot)$.  Generalizing a result of \cite{Bagyan} we obtain, reminiscent of Glivenko-Cantelli theory, the uniform convergence of $F_{Y_n|\Xi_n}(\cdot)$ to $F_{\mcalN_1(0,\sigma^2 V^2)}(\cdot)$ in the $p$th mean, for all $0 < p \le k$.  Further, we derive a Lipschitz continuity property of $F_{Y_n|\Xi_n}(\cdot) - F_{\mcalN_1(0,\sigma^2 V^2)}(\cdot)$.

In Section \ref{sec_examples_2}, we show that the additional requirements on $X_n$ in the main results in Sections \ref{sec_Lp_conv_cond_pdf} and \ref{sec_Lp_conv_cond_cdf} are satisfied by the examples studied in Section \ref{sec_examples}.  Further, we provide examples of random vectors $\Xi_n$ that satisfy the assumptions in our main results.  Finally, in Section \ref{sec_conclusions} we provide concluding remarks and discuss some open problems raised by our results.

\section{Some weak convergence properties of \texorpdfstring{$\boldsymbol{Y_n|\Xi_n}$}{YngivenXin}} 
\label{sec_gaussian_case}

Throughout this article the sequence $d_1,d_2,d_3,\ldots$ is monotonically increasing, and $d_n \to \infty$ as $n \to \infty$.  In denoting the dimensions by $d_n$ we are motivated by applications of random modulation in which $X_n$ and $\Xi_n$ are matrices; for instance, in compressed sensing \citep{Stern}, random modulation may involve $n \times n$ matrices (so that $d_n = n^2$) or $n \times n$ orthogonal matrices (in which case $d_n = \tfrac12 n(n-1)$).  

All vectors are column vectors, and all random vectors are continuous and have continuous density functions.  For a random entity $X$, we often write $\E_X$ to emphasize that the expectation is with respect to the marginal distribution of $X$.  Similarly, for any scalar random variable $Y$ and random entity $\Xi$, we denote by $\E_{Y|\Xi}$ and $\Var_{Y|\Xi}$ the mean and variance, respectively, with respect to the conditional distribution of $Y|\Xi$, and the conditional characteristic function of $Y|\Xi$ is $\varphi_{Y|\Xi}(t) = \E_{Y|\Xi} \exp(itY)$, $t \in \R$.

The probability distribution of $Y|\Xi$ is a \textit{random measure} \citep[Section 2]{Freedman}, and we use as the definition of \textit{weak convergence in probability} a characterization given by \citet[Lemma 2.2]{Diaconis_Freedman}:~For $n \in \N$ let $\mu_n$ be a random measure on $\R$ with (random) characteristic function $\widehat{\mu}_n$, and let $\mu_0$ be a deterministic measure on $\R$ with (deterministic) characteristic function $\widehat{\mu}_0$.  Then $\mu_n$ converges weakly in probability to $\mu_0$ as $n \to \infty$, denoted $\mu_n \cwip \mu_0$, if and only if $\widehat{\mu}_n(t) \cip \widehat{\mu}_0(t)$ for all $t \in \R$.

\subsection{Regularity conditions and weak convergence results for \texorpdfstring{$\boldsymbol{Y_n|\Xi_n}$}{YngivenXin}}
\label{subsec_regularity}

We assume throughout the article that the random vectors $\{X_n \in \R^{d_n}, n \ge 1\}$ satisfy the following regularity conditions: 

\begin{enumerate}[label=(C.\arabic*)]
\item \label{condition_1}
As $n \to \infty$, $\|X_n\|^2 \cip \sigma^2 > 0$. 
\item \label{condition_2}
Let $\wtilX_n$ be an independent copy of $X_n$.  Then $X_n'\wtilX_n \cip 0$ as $n \to \infty$.
\end{enumerate}

In stating these conditions, and throughout our work, we use the notation ``$X_n$'' in place of the more common scaling ``$X_n/\sqrt{d_n}$\,''.  With this notation duly noted, we remark that \ref{condition_1} and \ref{condition_2} are assumed widely in the literature.  \cite{Diaconis_Freedman} were first in stating \ref{condition_1} and \ref{condition_2} for the case in which $X_n$ has an empirical distribution; numerous authors (e.g., \cite{Bagyan}, \cite{Duembgen_etal}, \cite{Li_Yin}) formulated those assumptions subsequently for non-empirical distributions; and stronger versions of those conditions have been studied by other authors (e.g., \cite{Reeves}).  The conditions \ref{condition_1} and \ref{condition_2} have also appeared in the field of statistical physics \citep{Wee_Tatikonda}, where they are referred to as the ``thin-shell'' and ``zero overlap concentration'' assumptions, respectively.  

We write $\Xi \sim \mcalN_{d}(0,I_{d})$ to denote that a random vector $\Xi$ has a $d$-dimensional normal distribution with mean $0$ and covariance matrix $I_{d}$, the identity matrix of order $d$.  We also use the notation $\i = \sqrt{-1}$, and we often write $\E_X$ to emphasize that an expectation is being taken with respect to the distribution of a given random entity $X$.  

Let $k$ and $l$ be fixed positive integers, and let $X_{n,1},\ldots,X_{n,k} \in \R^{d_n}$ be mutually independent, each satisfying \ref{condition_1} and \ref{condition_2}.  Also let $\Xi_{n,1},\ldots,\Xi_{n,l}$ be mutually independent copies of $\Xi_n \sim \mcalN_{d_n}(0,I_{d_n})$, with $\{\Xi_{n,1},\ldots,\Xi_{n,l}\}$ and $\{X_{n,1},\ldots,X_{n,k}\}$ also are independent.  This situation arises in practice when, given a random sample $X_{n,1},\ldots,X_{n,k}$ from $X_n$, we simulate a random sample $\Xi_{n,1},\ldots,\Xi_{n,l}$ from $\Xi_n$ and then seek to use the family of projections $Y_{n;j,r} = \Xi_{n,j}'X_{n,r}$, $j=1,\ldots,l$, $r=1,\ldots,k$, to perform inference for the population represented by the conditional distribution of $\Xi_n'X_n|\Xi_n$.  

Defining the $l \times k$ matrix $\mcalY_n = (Y_{n;j,r})$, we now provide the asymptotic conditional distribution of $\mcalY_n$, given $(\Xi_{n,1},\ldots,\Xi_{n,l})$, as $n \to \infty$.  

\begin{theorem}
\label{thm_resampling_Xi_n_X_n}
Let $X_{n,1},\ldots,X_{n,k} \in \R^{d_n}$, $n \ge 1$, be mutually independent copies of $X_n$, where $X_n$ satisfies \ref{condition_1} and \ref{condition_2}.  Let $\Xi_{n,1},\ldots,\Xi_{n,l} \in \R^{d_n}$ be mutually independent $\mcalN_{d_n}(0,I_{d_n})$ vectors, all independent of $(X_{n,1},\ldots,X_{n,k})$.  Then $\mcalY_n|(\Xi_{n,1},\ldots,\Xi_{n,l}) \cwip \mcalZ$ as $n \to \infty$, where $\mcalZ = (Z_{j,r})$ is an $l \times k$ random matrix whose entries $Z_{j,r}$, $j=1,\ldots,l$, $r=1,\ldots,k$, are mutually independent and $\mcalN_1(0,\sigma^2)$--distributed.  
\end{theorem}

For the case in which $k = 1$, Theorem \ref{thm_resampling_Xi_n_X_n} reduces to the following result of \citet[Corollary 2.2]{Duembgen_etal}.  

\begin{corollary} {\rm{\citep{Duembgen_etal}}}
\label{cor_modulation_normal_distn}
For each $n \in \N$, suppose that $X_n \in \R^{d_n}$ satisfies \ref{condition_1} and \ref{condition_2}.  Let the random vectors $\Xi_{n,1},\ldots,\Xi_{n,l} \in \R^{d_n}$ be mutually independent, $\mcalN_{d_n}(0,I_{d_n})$--distributed, and independent of $X_n$; and define $\mcalY_n = (\Xi_{n,1}'X_n,\ldots,\Xi_{n,l}'X_n)'$.  Then $\mcalY_n|(\Xi_{n,1},\ldots,\Xi_{n,l}) \cwip \mcalN_l(0,\sigma^2 I_l)$ as $n \to \infty$.  
\end{corollary}

\begin{remark}
\label{rem_caution}
{\rm 
(i) \citet[Lemma 4.1]{Duembgen_etal} also proved the following converse to Corollary \ref{cor_modulation_normal_distn}, the proof of which can be readily adapted to our setting: {Suppose that $\{X_n \in \R^{d_n}, n \ge 1\}$ and $\{\Xi_n \in \R^{d_n}, n \ge 1\}$ are independent, and let $Y_n = \Xi_n'X_n$.  If $Y_n|\Xi_n \cwip \mcalN_1(0,\sigma^2)$ as $n \to \infty$ then \ref{condition_1} and \ref{condition_2} hold.}  

(ii) By Corollary \ref{cor_modulation_normal_distn} $Y_n|\Xi_n \cwip \mcalN_1(0,\sigma^2)$, which does not depend on $\Xi_n$, so the corresponding unconditional distribution of $Y_n$ also converges similarly to $\mcalN_1(0,\sigma^2)$.  This property, in which the limiting conditional distribution of $Y_n|\Xi_n$ does not depend on $\Xi_n$, appears repeatedly in the sequel.  

(iii) As noted by \citet[p.~94]{Duembgen_etal}, results such as Corollary \ref{cor_modulation_normal_distn} caution us to be wary of presuming, on the basis of moderately many low-dimensional projections, that a high-dimensional data set is normally distributed.  

(iv) In much of the literature, $X_n$ is projected along uniformly distributed directions.  To recover this case from our results, one sets $\Xi_n = \sqrt{d_n}\Theta_n$ where $\Theta_n$ is uniformly distributed on $\mcalS^{d_n - 1}$, the hypersphere centered at the origin and of radius $1$.  Then $\Xi_n'X_n = \sqrt{d_n} \Theta_n'X_n \eqdist \sqrt{d_n} \Theta_{n,1} \|X\|$, $\Theta_{n,1}$ being the first component of $\Theta_n$, and the proof of Corollary \ref{cor_modulation_normal_distn} carries over, using the fact that the distribution of $\sqrt{d_n}\Theta_{n,1}$ converges uniformly to a standard normal distribution.  
}\end{remark}

\begin{remark}
\label{rem_non_Gaussian_Xin}
{\rm There is the issue of whether Theorem \ref{thm_resampling_Xi_n_X_n} can be extended to the case in which $\Xi_n$ has a non-Gaussian distribution.  Consider, for simplicity, the case in which $l=k=1$; then it will be shown that 
\begin{equation}
\label{eq_Var_to_0}
\Var_{\Xi_n}\big(\varphi_{Y_n|\Xi_n}(t)\big) := \E_{\Xi_n} \big|\varphi_{Y_n|\Xi_n}(t)\big|^2 - \big|\E_{\Xi_n} \varphi_{Y_n|\Xi_n}(t)\big|^2 \to 0
\end{equation}
as $n \to \infty$, so this raises the issue of whether \eqref{eq_Var_to_0} holds for non-Gaussian $\Xi_n$.  

Suppose that $\Xi_n$ has a spherically symmetric stable distribution with \textit{index of stability} $\alpha \in (0,2)$ \citep{Zolotarev} and characteristic function $\E \exp(\i u'\Xi_n) = \exp(- \|u\|^\alpha)$, $u \in \R^{d_n}$.   Then it will be shown in Subsection \ref{subsec_regularity_proofs}, starting at \eqref{eq_E_phi_n_stable}, that 
\begin{equation}
\label{eq_Var_symm_stable}
\lim_{n \to \infty} \Var_{\Xi_n}\big(\varphi_{Y_n|\Xi_n}(t)\big) = \exp(- 2^{\alpha/2} \sigma^\alpha |t|^\alpha) - \exp(- 2 \sigma^\alpha |t|^\alpha),
\end{equation}
which is positive for $t \neq 0$, so the proof of Theorem \ref{thm_resampling_Xi_n_X_n} does not apply in this case.  

It is noticeable that the distribution of $\Xi_n$ in this counterexample is \textit{spherically symmetric}, \textit{i.e.}, the characteristic function $\E \exp(\i u'\Xi_n)$ is a function of $\|u\|$.  This also raises the issue of the extent to which \eqref{eq_Var_to_0} is characteristic of the normal distribution, and indeed we show that, subject to \ref{condition_1} and \ref{condition_2}, the property \eqref{eq_Var_to_0} characterizes the normality of $\Xi_n$ within the class of spherically symmetric distributions.  
}\end{remark}

\begin{theorem}
\label{thm_modulation_charac_2}
Suppose that $\{X_n \in \R^{d_n}, n \ge 1\}$ satisfy \ref{condition_1} and \ref{condition_2}, and let $\{\Xi_n \in \R^{d_n}, n \ge 1\}$ be mutually independent of $\{X_n, n \ge 1\}$.  Also suppose that $\E \exp(\i u'\Xi_n) = \psi_0(\|u\|^2)$, $u \in \R^{d_n}$, for some function $\psi_0: [0,\infty) \to \R$, 
and define $Y_n = \Xi_n'X_n$, $n \ge 1$.  Then $\Xi_n \sim \mcalN_{d_n}(0,\sigma_0^2 I_{d_n})$ for some $\sigma_0$ if and only if \eqref{eq_Var_to_0} holds.  
\end{theorem}

\begin{remark}
\label{rem_modulation_charac_2}
{\rm
Suppose that $Z_1$ and $Z_2$ are independent copies of a random variable $Z$, and $Z \eqdist 2^{-1/2} (Z_1+Z_2)$.  By a theorem of \cite{Polya} (cf., \citet[Theorem 1.9.5]{Bogachev}), $Z \sim \mcalN_1(0,\sigma_0^2)$ for some $\sigma_0$.   

There is an extensive literature that proves P\'olya's theorem to be ``stable,'' \textit{i.e.}, if $Z$ and $2^{-1/2} (Z_1+Z_2)$ are ``close in distribution'' according to various measures of closeness, then $Z$ is close in distribution to $\mcalN_1(0,\sigma_0^2)$ \citep{YY2}.  Extensions of Corollary \ref{cor_modulation_normal_distn} and Theorem \ref{thm_modulation_charac_2} can be obtained from such stability results, and we leave such details to interested readers.  
}\end{remark}

\subsection{Proofs}
\label{subsec_regularity_proofs}

\noindent
\textit{Proof of Theorem \ref{thm_resampling_Xi_n_X_n}}:  Let $U = (U_{j,r})$, a constant $l \times k$ real matrix, and define 
$$
Z_n = \tr(U'\mcalY_n) \equiv \sum_{j=1}^l \sum_{r=1}^k u_{j,r} Y_{n;j,r} = \sum_{j=1}^l \sum_{r=1}^k u_{j,r} \Xi_{n,j}'X_{n,r}.
$$
By the mutual independence of $\Xi_{n,1},\ldots,\Xi_{n,l}$, their independence from $X_{n,1},\ldots,X_{n,k}$, and Fubini's theorem, we obtain, for any $t \in \R$, 
\begin{align*}
\E_{(\Xi_{n,1},\ldots,\Xi_{n,l})} \varphi_{Z_n|(\Xi_{n,1},\ldots,\Xi_{n,l})}(t) 
&= \E_{(\Xi_{n,1},\ldots,\Xi_{n,l})} \E_{(X_{n,1},\ldots,X_{n,k})} \exp\Big(\i t \sum_{j=1}^l \sum_{r=1}^k u_{j,r} \Xi_{n,j}'X_{n,r}\Big) \\
&= \E_{(X_{n,1},\ldots,X_{n,k})} \E_{(\Xi_{n,1},\ldots,\Xi_{n,l})} \exp\Big(\i t \sum_{j=1}^l \sum_{r=1}^k u_{j,r} \Xi_{n,j}'X_{n,r}\Big).
\end{align*}
Since $\Xi_{n,j} \sim \mcalN_{d_n}(0,I_{d_n})$, $j=1,\ldots,l$, then it follows that 
\begin{align*}
\E_{(\Xi_{n,1},\ldots,\Xi_{n,l})} \varphi_{Z_n|(\Xi_{n,1},\ldots,\Xi_{n,l})}(t) 
&= \E_{(X_{n,1},\ldots,X_{n,k})} \prod_{j=1}^l \E_{\Xi_{n,j}} \exp\Big(\i t \Xi_{n,j}' \sum_{r=1}^k u_{j,r} X_{n,r}\Big) \\
&= \E_{(X_{n,1},\ldots,X_{n,k})} \prod_{j=1}^l \exp\Big(-\tfrac12 t^2 \Big\|\sum_{r=1}^k u_{j,r} X_{n,r}\Big\|^2\Big).
\end{align*}
Denoting Kronecker's delta by $\delta_{j,r}$, it follows from \ref{condition_1} and \ref{condition_2} that 
\begin{equation}
\label{eq_norm_sum_ujrXnr}
\Big\|\sum_{r=1}^k u_{j,r} X_{n,r}\Big\|^2 = \sum_{r_1=1}^k \sum_{r_2=1}^k u_{j,r_1} u_{j,r_2} X_{n,r_1}'X_{n,r_2} \cip \sigma^2 \sum_{r=1}^k u_{j,r}^2.
\end{equation}
It follows by the continuity of the exponential function and the Continuous Mapping Theorem \cite[p.~254, Theorem 1]{Chow_Teicher} that, as $n \to \infty$, 
\begin{equation}
\label{E_phi_n_limit}
\E_{(\Xi_{n,1},\ldots,\Xi_{n,l})} \varphi_{Z_n|(\Xi_{n,1},\ldots,\Xi_{n,l})}(t) \to \prod_{j=1}^l \exp\Big(-\tfrac12 t^2 \sigma^2 \sum_{r=1}^k u_{j,r}^2\Big) 
\equiv \exp(-\tfrac12 t^2 \sigma^2 \tr U'U).
\end{equation}

Let $\wtilX_{n,1},\ldots,\wtilX_{n,k}$ be mutually independent copies of $X_{n,1},\ldots,X_{n,k}$.  Then 
\begin{align*}
\E&{}_{(\Xi_{n,1},\ldots,\Xi_{n,l})} \big|\varphi_{Z_n|(\Xi_{n,1},\ldots,\Xi_{n,l})}(t)\big|^2 \\
&\ \ \ = \E_{(\Xi_{n,1},\ldots,\Xi_{n,l})} \left[\varphi_{Z_n|(\Xi_{n,1},\ldots,\Xi_{n,l})}(t) \, \overline{\varphi_{Z_n|(\Xi_{n,1},\ldots,\Xi_{n,l})}(t)}\right] \nonumber \\
&\ \ \ = \E_{(\Xi_{n,1},\ldots,\Xi_{n,l})} \bigg[\E_{X_{n,1},\ldots,X_{n,k}|(\Xi_{n,1},\ldots,\Xi_{n,l})} \exp\Big(\i t \sum_{j=1}^l \sum_{r=1}^k u_{j,r} \Xi_{n,j}'X_{n,r}\Big) \\
&\ \ \ \qquad\qquad\qquad\ \cdot \E_{\wtilX_{n,1},\ldots,\wtilX_{n,k}|(\Xi_{n,1},\ldots,\Xi_{n,l})} \exp\Big(-\i t \sum_{j=1}^l \sum_{r=1}^k u_{j,r} \Xi_{n,j}'\wtilX_{n,r}\Big)\bigg].
\end{align*}
Reversing the order of the expectations, which is justified by Fubini's theorem, and applying the mutual independence of $\Xi_{n,1},\ldots,\Xi_{n,l}$, we obtain 
\begin{align}
\label{eq_Y_n_Xi_n}
&\E{}_{(\Xi_{n,1},\ldots,\Xi_{n,l})} \big|\varphi_{Z_n|(\Xi_{n,1},\ldots,\Xi_{n,l})}(t)\big|^2 \nonumber \\
&\ = \E_{X_{n,1},\ldots,X_{n,k}} \E_{\wtilX_{n,1},\ldots,\wtilX_{n,k}} \prod_{j=1}^l \E_{\Xi_{n,j}} \exp\Big(\i t \Xi_{n,j}'\sum_{r=1}^k u_{j,r} (X_{n,r} - \wtilX_{n,r})\Big) \nonumber \\
&\ = \E_{X_{n,1},\ldots,X_{n,k}} \E_{\wtilX_{n,1},\ldots,\wtilX_{n,k}} \prod_{j=1}^l \exp\Big(- \tfrac12 t^2 \Big\|\sum_{r=1}^k u_{j,r} (X_{n,r} - \wtilX_{n,r})\Big\|^2\Big).
\end{align}

It is straightforward that 
\begin{equation}
\label{eq_norm_sum}
\Big\|\sum_{r=1}^k u_{j,r} (X_{n,r} - \wtilX_{n,r})\Big\|^2 = \sum_{r_1=1}^k \sum_{r_2=1}^k u_{j,r_1} u_{j,r_2} (X_{n,r_1} - \wtilX_{n,r_1})'(X_{n,r_2} - \wtilX_{n,r_2}).
\end{equation}
By \ref{condition_1}, \ref{condition_2}, the mutual independence of $\{X_{n,1},\ldots,X_{n,k}\}$ and $\{\wtilX_{n,1},\ldots,\wtilX_{n,k}\}$, and Slutsky's theorem, it follows that  
\begin{multline}
\label{eq_inner_X_Xtilde}
(X_{n,r_1} - \wtilX_{n,r_1})'(X_{n,r_2} - \wtilX_{n,r_2}) \\
= X_{n,r_1}'X_{n,r_2} - X_{n,r_1}' \wtilX_{n,r_2} - \wtilX_{n,r_1}' X_{n,r_2} + \wtilX_{n,r_1}' \wtilX_{n,r_2} \cip 2 \delta_{r_1,r_2} \sigma^2
\end{multline}
as $n \to \infty$.  Applying \eqref{eq_inner_X_Xtilde} to \eqref{eq_norm_sum}, we obtain 
$$
\Big\|\sum_{r=1}^k u_{j,r} (X_{n,r} - \wtilX_{n,r})\Big\|^2 \cip 2\sigma^2 \sum_{r=1}^k u_{j,r}^2,
$$
and it follows by the Continuous Mapping Theorem that 
\begin{equation}
\label{E_phi_n_sq_limit}
\E_{(\Xi_{n,1},\ldots,\Xi_{n,l})} \big|\varphi_{Z_n|(\Xi_{n,1},\ldots,\Xi_{n,l})}(t)\big|^2 \to \prod_{j=1}^l \exp\Big(- t^2 \sigma^2 \sum_{r=1}^k u_{j,r}^2\Big) = \exp(- t^2 \sigma^2 \tr U'U).
\end{equation}

Next, for $\varepsilon > 0$, it follows by Chebyshev's inequality that 
\begin{align*}
\P&\big(\big|\varphi_{Z_n|(\Xi_{n,1},\ldots,\Xi_{n,l})}(t) - \exp(- \tfrac12 t^2 \sigma^2 \tr U'U)\big| > \varepsilon\big) \\
&\le \varepsilon^{-2} \E_{(\Xi_{n,1},\ldots,\Xi_{n,l})} \big|\varphi_{Z_n|(\Xi_{n,1},\ldots,\Xi_{n,l})}(t) - \exp(- \tfrac12 t^2 \sigma^2 \tr U'U)\big|^2 \\
&\equiv \varepsilon^{-2} \E_{(\Xi_{n,1},\ldots,\Xi_{n,l})} \Big[\big|\varphi_{Z_n|(\Xi_{n,1},\ldots,\Xi_{n,l})}(t)\big|^2 - \exp(- t^2 \sigma^2 \tr U'U) \\
& \qquad\qquad\qquad\qquad - \big(\varphi_{Z_n|(\Xi_{n,1},\ldots,\Xi_{n,l})}(t) - \exp(- \tfrac12 t^2 \sigma^2 \tr U'U)\big) \exp(- \tfrac12 t^2 \sigma^2 \tr U'U) \\
& \qquad\qquad\qquad\qquad\ -  \big(\overline{\varphi_{Z_n|(\Xi_{n,1},\ldots,\Xi_{n,l})}(t)} - \exp(- \tfrac12 t^2 \sigma^2 \tr U'U)\big) \exp(- \tfrac12 t^2 \sigma^2 \tr U'U)\Big].
\end{align*}
Applying the triangle inequality, and the inequality $\exp(- \tfrac12 t^2 \sigma^2 \tr U'U) \le 1$ for all $t$ and $U$, we obtain 
\begin{align}
\label{eq_chebyshev}
\P\big(\big|\varphi&{}_{Z_n|(\Xi_{n,1},\ldots,\Xi_{n,l})}(t) - \exp(- \tfrac12 t^2 \sigma^2 \tr U'U)\big| > \varepsilon\big) \nonumber \\
&\le \varepsilon^{-2} \Big[\big|\E_{(\Xi_{n,1},\ldots,\Xi_{n,l})} \big|\varphi_{Z_n|(\Xi_{n,1},\ldots,\Xi_{n,l})}(t)\big|^2 - \exp(- t^2 \sigma^2 \tr U'U)\big| \nonumber \\
& \qquad\quad + \big|\E_{(\Xi_{n,1},\ldots,\Xi_{n,l})} \varphi_{Z_n|(\Xi_{n,1},\ldots,\Xi_{n,l})}(t) - \exp(- \tfrac12 t^2 \sigma^2 \tr U'U)\big| \nonumber \\
& \qquad\qquad\ + \big|\E_{(\Xi_{n,1},\ldots,\Xi_{n,l})} \overline{\varphi_{Z_n|(\Xi_{n,1},\ldots,\Xi_{n,l})}(t)} - \exp(- \tfrac12 t^2 \sigma^2 \tr U'U)\big|\Big].
\end{align}
By \eqref{E_phi_n_limit} and \eqref{E_phi_n_sq_limit}, each of the three terms on the right-hand side of \eqref{eq_chebyshev} converges to $0$ as $n \to \infty$.  Since $\epsilon$ was chosen arbitrarily then it follows that, for all $t$ and $U$, 
\begin{equation}
\label{eq_Zn_given_Xin}
\varphi_{Z_n|(\Xi_{n,1},\ldots,\Xi_{n,l})}(t) \cip \exp(- \tfrac12 t^2 \sigma^2 \tr U'U),
\end{equation}
the characteristic function of the $\mcalN_1(0,\sigma^2 \tr U'U)$ distribution.  
Applying the characterization of weak convergence in probability given 
at Subsection \ref{subsec_regularity}, \textit{supra}, 
we obtain $Z_n|(\Xi_{n,1},\ldots,\Xi_{n,l}) \cwip \mcalN_1(0,\sigma^2 \tr U'U)$.  Finally, since $U$ was chosen arbitrarily then it follows by the Cram\'er-Wold device that $\mcalY_n|(\Xi_{n,1},\ldots,\Xi_{n,l}) \cwip \mcalZ$.
$\qed$

\medskip

\noindent
\textit{Proof of Equation \eqref{eq_Var_symm_stable}}:  
By an interchange of expectations, which is valid by an application of Fubini's theorem, and using the independence of $X_n$ and $\Xi_n$, we obtain 
\begin{equation}
\label{eq_E_phi_n_stable}
\E_{\Xi_n} \varphi_{Y_n|\Xi_n}(t) = \E_{X_n} \E_{\Xi_n} \exp(\i t\Xi_n'X_n) = \E \exp(- |t|^\alpha \|X_n\|^\alpha),
\end{equation}
$t \in \R$.  By \ref{condition_1}, the continuity of the exponential function, and the Continuous Mapping Theorem, it follows from \eqref{eq_E_phi_n_stable} that, for all $t \in \R$, 
\begin{equation}
\label{eq_E_phi_n_sq_limit_stable_2}
\lim_{n \to \infty} |\E_{\Xi_n} \varphi_{Y_n|\Xi_n}(t)|^2 =  [\exp(- \sigma^\alpha |t|^\alpha)]^2 = \exp(- 2 \sigma^\alpha |t|^\alpha).
\end{equation}
Let $\wtilX_n$ be an independent copy of $X_n$.  Similar to \eqref{eq_Y_n_Xi_n}, we obtain 
\begin{equation}
\label{eq_Y_n_Xi_n_stable}
\E_{\Xi_n} \big|\varphi_{Y_n|\Xi_n}(t)\big|^2 = \E_{X_n,\wtilX_n} \exp\big(- |t|^\alpha \|X_n - \wtilX_n\|^\alpha\big).
\end{equation}
By \eqref{eq_inner_X_Xtilde}, 
$
\|X_n - \wtilX_n\|^\alpha \equiv (\|X_n - \wtilX_n\|^2)^{\alpha/2} \cip 2^{\alpha/2} \sigma^\alpha,
$ 
and by Slutsky's theorem, 
$$
\E_{X_n,\wtilX_n} \exp\big(- |t|^\alpha \|X_n - \wtilX_n\|^\alpha\big) \to \exp(- 2^{\alpha/2} \sigma^\alpha |t|^\alpha)
$$
as $n \to \infty$.  Applying \eqref{eq_Y_n_Xi_n_stable} we obtain, for all $t \in \R$, 
\begin{equation}
\label{eq_E_phi_n_sq_limit_stable}
\lim_{n \to \infty} \E_{\Xi_n} \big|\varphi_{Y_n|\Xi_n}(t)\big|^2 = \exp(- 2^{\alpha/2} \sigma^\alpha |t|^\alpha),
\end{equation}
and by combining \eqref{eq_E_phi_n_sq_limit_stable} and \eqref{eq_E_phi_n_sq_limit_stable_2} we obtain \eqref{eq_Var_symm_stable}.  
$\qed$

\bigskip

\noindent
\textit{Proof of Theorem \ref{thm_modulation_charac_2}}:  
Suppose that $\Xi_n \sim \mcalN_{d_n}(0,\sigma_0^2 I_{d_n})$ for some $\sigma_0$.  By \eqref{E_phi_n_limit}, $|\E_{\Xi_n} \varphi_{Y_n|\Xi_n}(t)|^2 \to \exp(- t^2 \sigma_0^2)$.  Also, by \eqref{E_phi_n_sq_limit}, $\E_{\Xi_n} |\varphi_{Y_n|\Xi_n}(t)|^2 \to \exp(- t^2 \sigma_0^2)$.  
Therefore $\Var_{\Xi_n}\big(\varphi_{Y_n|\Xi_n}(t)\big) \to 0$ as $n \to \infty$.  

Conversely suppose that, for all $t \in \R$, $\Var_{\Xi_n}\big(\varphi_{Y_n|\Xi_n}(t)\big) \to 0$ as $n \to \infty$.  Then  
\begin{align*}
\E_{\Xi_n} \varphi_{Y_n|\Xi_n}(t) &= \E_{\Xi_n} \E_{Y_n|\Xi_n} \exp(\i tY_n) \\
&= \E_{X_n} \E_{\Xi_n} \exp(\i tX_n'\Xi_n) = \E_{X_n} \psi_0(t^2 \|X_n\|^2).
\end{align*}
Since $\psi_0(t^2)$ is a characteristic function then it is continuous.  By \ref{condition_1}, $\|X_n\|^2 \cip \sigma^2$, so by the Continuous Mapping Theorem, $\psi_0(t^2 \|X_n\|^2) \cip \psi_0(t^2 \sigma^2)$, $t \in \R$; therefore 
\begin{equation}
\label{eq_Y_n_Xi_n_15}
\lim_{n \to \infty} \E_{\Xi_n} \varphi_{Y_n|\Xi_n}(t) = \psi_0(\sigma^2 t^2).
\end{equation}

Let $\wtilX_n$ be an independent copy of $X_n$.  As $\varphi_{Y_n|\Xi_n}(\cdot)$ is a characteristic function then it is bounded, so by applying Fubini's theorem to interchange expectations we obtain 
\begin{align*}
\E_{\Xi_n} \big|\varphi_{Y_n|\Xi_n}(t)\big|^2 &= \E_{\Xi_n} \big[\varphi_{Y_n|\Xi_n}(t) \, \overline{\varphi_{Y_n|\Xi_n}(t)}\,\big] \nonumber \\
&= \E_{\Xi_n} \E_{X_n|\Xi_n} \exp(\i t \Xi_n'X_n) \cdot \E_{\wtilX_n|\Xi_n} \exp(-\i t \Xi_n'\wtilX_n) \\
&= \E_{\Xi_n} \E_{X_n,\wtilX_n} \exp\big(\i t \Xi_n(X_n - \wtilX_n)\big),
\end{align*}
where the latter equality follows from the law of iterated expectations.  Again interchanging expectations, which is justified by Fubini's theorem, we obtain 
\begin{align*}
\E_{\Xi_n} \big|\varphi_{Y_n|\Xi_n}(t)\big|^2 &= \E_{X_n,\wtilX_n} \E_{\Xi_n} \exp\big(\i t \Xi_n(X_n - \wtilX_n)\big) 
= \E_{X_n,\wtilX_n} \psi_0(t^2 \|X_n - \wtilX_n\|^2\big).
\end{align*}
By \ref{condition_1}, \ref{condition_2}, and Slutsky's theorem, $\|X_n-\wtilX_n\|^2 \cip 2\sigma^2$ as $n \to \infty$.  Since $\psi_0(\cdot)$ is continuous then, by the Continuous Mapping Theorem, 
\begin{equation}
\label{eq_Y_n_Xi_n_2}
\lim_{n \to \infty} \E_{\Xi_n} \big|\varphi_{Y_n|\Xi_n}(t)\big|^2 = \psi_0(2\sigma^2 t^2).
\end{equation}
Combining \eqref{eq_Y_n_Xi_n_15} and \eqref{eq_Y_n_Xi_n_2}, we obtain 
\begin{align*}
\psi_0(2\sigma^2 t^2) - [\psi_0(\sigma^2 t^2)]^2 &= \lim_{n \to \infty} \big[\E_{\Xi_n} \big|\varphi_{Y_n|\Xi_n}(t)\big|^2 - \big|\E_{\Xi_n} \varphi_{Y_n|\Xi_n}(t)\big|^2 \big] \\
&= \lim_{n \to \infty} \Var_{\Xi_n}\big(\varphi_{Y_n|\Xi_n}(t)\big) = 0.
\end{align*}
Therefore we obtain the functional equation, $\psi_0(2\sigma^2 t^2) = [\psi_0(\sigma^2 t^2)]^2$, equivalently, 
\begin{equation}
\label{eq_Polya_eq}
\psi_0(t^2) = [\psi_0(2^{-1} t^2)]^2, \qquad t \in \R.
\end{equation}

Denote by $Z$ a random variable with characteristic function $\psi_0(t^2)$, and let $Z_1$ and $Z_2$ be mutually independent random variables that have the same distribution as $Z$.  Then the right-hand side of \eqref{eq_Polya_eq} is the characteristic function of $2^{-1/2} (Z_1 + Z_2)$, so \eqref{eq_Polya_eq} is equivalent to the equality in distribution, 
\begin{equation}
\label{eq_decomp_device}
Z \eqdist 2^{-1/2} (Z_1 + Z_2).
\end{equation}
By P\'olya's theorem, \eqref{eq_decomp_device} implies that $Z \sim \mcalN_1(0,\sigma_0^2)$ for some $\sigma_0$, so $\psi_0(t^2) = \exp(-\tfrac12 \sigma_0^2 t^2)$, $t \in \R$.  Therefore $\E_{\Xi_n} \exp(\i u'\Xi_n) = \psi_0(\|u\|^2) = \exp(-\tfrac12 \sigma_0^2 \|u\|^2)$, $u \in \R^{d_n}$, hence $\Xi_n \sim \mcalN_{d_n}(0,\sigma_0^2 I_{d_n})$.  
$\qed$

\smallskip

\begin{remark}
\label{rem_C_1_prime_sec_2}
{\rm
We are motivated by a comment of the reviewer to determine the effects of replacing \ref{condition_1} by 
\begin{enumerate}[label=(C.\arabic*)$^\prime$]
\item \label{condition_1_prime}
As $n \to \infty$, $\|X_n\|^2 \cip R^2$, where the random variable $R \neq 0$, almost surely.
\end{enumerate}
To address this issue, we consider 
the proof of Theorem \ref{thm_resampling_Xi_n_X_n} and modify the calculations accordingly.  Under \ref{condition_1_prime}, \eqref{eq_norm_sum_ujrXnr} becomes 
$$
\Big\|\sum_{r=1}^k u_{j,r} X_{n,r}\Big\|^2 = \sum_{r_1=1}^k \sum_{r_2=1}^k u_{j,r_1} u_{j,r_2} X_{n,r_1}'X_{n,r_2} \cip \sum_{r=1}^k u_{j,r}^2 R_r^2,
$$
where $R_1,\ldots,R_k$ are independent copies of $R$.  Define $\mcalR_k = \diag(R_1^2,\ldots,R_k^2)$, then \eqref{E_phi_n_limit} becomes 
\begin{eqnarray*}
\E_{(\Xi_{n,1},\ldots,\Xi_{n,l})} \varphi_{Z_n|(\Xi_{n,1},\ldots,\Xi_{n,l})}(t) & \to & \E_{R_1,\ldots,R_k} \prod_{j=1}^l \exp\Big(-\tfrac12 t^2 \sum_{r=1}^k u_{j,r}^2 R_r^2\Big) \\
& = & \E_{R_1,\ldots,R_k} \exp(-\tfrac12 t^2 \tr U'U\mcalR_k),
\end{eqnarray*}
and, similarly, \eqref{E_phi_n_sq_limit} becomes 
$$
\E_{(\Xi_{n,1},\ldots,\Xi_{n,l})} \big|\varphi_{Z_n|(\Xi_{n,1},\ldots,\Xi_{n,l})}(t)\big|^2 \to \E \exp(- t^2 \tr U'U\mcalR_k).
$$
Next, \eqref{eq_Zn_given_Xin} becomes 
\begin{equation}
\label{eq_Zn_given_Xin_R}
\varphi_{Z_n|(\Xi_{n,1},\ldots,\Xi_{n,l})}(t) \cip \E \exp(- \tfrac12 t^2 \tr U'U\mcalR_k),
\end{equation}
which is the characteristic function of the $l \times k$ random matrix $\mcalZ = (R_r^2 Z_{j,r})$ where $R_1,\ldots,R_k$ are as before; each $Z_{j,r} \sim \mcalN_1(0,1)$; and $R_1,\ldots,R_k,Z_{1,1},\ldots,Z_{l,k}$ are mutually independent.  Therefore $\mcalY_n|(\Xi_{n,1},\ldots,\Xi_{n,l}) \cwip \mcalZ$ as $n \to \infty$.  

Next, we generalize \eqref{eq_Var_symm_stable} under the assumption of \ref{condition_1_prime}.  Similar to \eqref{eq_E_phi_n_stable} and \eqref{eq_E_phi_n_sq_limit_stable_2} there holds, for all $t \in \R$,  
$$
\lim_{n \to \infty} \E_{\Xi_n} \varphi_{Y_n|\Xi_n}(t) =  \E \exp(- R^\alpha |t|^\alpha).
$$
Therefore, with $R_1, R_2$ independent and distributed as $R$, 
\begin{equation}
\label{eq_sq_E_phi_n_limit_stable_R}
\lim_{n \to \infty} |\E_{\Xi_n} \varphi_{Y_n|\Xi_n}(t)|^2 = \Big|\lim_{n \to \infty} \E_{\Xi_n} \varphi_{Y_n|\Xi_n}(t)\Big|^2 = \E \exp\big( \! - (R_1^\alpha + R_2^\alpha) |t|^\alpha\big),
\end{equation}
By modifying the arguments that led to \eqref{eq_Y_n_Xi_n_stable}-\eqref{eq_E_phi_n_sq_limit_stable}, and then applying Slutsky's theorem, we obtain 
\begin{equation}
\label{eq_E_sq_phi_n_limit_stable_R}
\E_{\Xi_n} \big|\varphi_{Y_n|\Xi_n}(t)\big|^2 \to \E \exp\big(\!- (R_1^2 + R_2^2)^{\alpha/2} |t|^\alpha\big)
\end{equation}
as $n \to \infty$.  Combining \eqref{eq_sq_E_phi_n_limit_stable_R} and \eqref{eq_E_sq_phi_n_limit_stable_R}, we obtain 
\begin{equation}
\label{eq_Var_phi_n_sq_limit_stable_R}
\lim_{n \to \infty} \Var_{\Xi_n} \big(\varphi_{Y_n|\Xi_n}(t)\big) = \E \big[\exp\big(\!- (R_1^2 + R_2^2)^{\alpha/2} |t|^\alpha\big) - \exp\big(\!- (R_1^\alpha + R_2^\alpha) |t|^\alpha\big)\big].
\end{equation}
It simple to show that 
$(r_1^2 + r_2^2)^{\alpha/2} < r_1^\alpha + r_2^\alpha$ 
for all $r_1, r_2 > 0$ and $0 < \alpha < 2$, so it follows that the right-hand side of \eqref{eq_Var_phi_n_sq_limit_stable_R} is positive for all $t \neq 0$.  

Next, we extend Theorem \ref{thm_modulation_charac_2} under the assumption \ref{condition_1_prime}.  Suppose that $\Xi_n \sim \mcalN_{d_n}(0,\sigma_0^2 I_{d_n})$, $\sigma_0 > 0$; then by \eqref{eq_Zn_given_Xin_R} with $k=l=1$, $\varphi_{Y_n|\Xi_n}(t)\big) \cip \E \exp(-\tfrac12 t^2 R^2)$
as $n \to \infty$, and therefore $\Var_{\Xi_n} \big(\varphi_{Y_n|\Xi_n}(t)\big) \to 0$.  
Conversely, suppose that $\Xi_n$ is spherically symmetric with $\E \exp(\i u'\Xi_n) = \psi_0(\|u\|^2)$, $u \in \R^{d_n}$; proceeding similarly to the derivations of \eqref{eq_Y_n_Xi_n_15} and \eqref{eq_Y_n_Xi_n_2}, we obtain 
\begin{align}
\label{eq_int_eq_char_R}
\Var_{\Xi_n} \big(\varphi_{Y_n|\Xi_n}(t)\big) &= \E \psi_0\big(\|X_n-\wtilX_n\|^2 t^2\big) - \big[\E\psi_0(\|X_n\|^2 t^2)\big]^2 \nonumber \\
&\to \E \psi_0\big((R_1^2 + R_2^2) t^2\big) - \big[\E \psi_0(R^2 t^2)\big]^2.
\end{align}
In an article under preparation \citep{Bagyan_Richards_26}, we have generalized P\'olya's characterization of the normal distribution by showing that for certain characteristic functions $\psi_0(t^2)$ that are analytic at $t=0$, \eqref{eq_int_eq_char_R} is identically $0$ if and only if $\psi_0(t^2) = \exp(-\tfrac12 \sigma_0^2 t^2)$ for some $\sigma_0 > 0$; this holds, for example, if $R^2$ has finitely many possible values.  For such random variables $R$, Theorem \ref{thm_modulation_charac_2} remains valid under \ref{condition_1_prime}.  
}\end{remark}

\section{Examples of distributions satisfying \ref{condition_1} and \ref{condition_2}}
\label{sec_examples}

We now provide examples to illustrate the breadth of the class of distributions that satisfy  \ref{condition_1} and \ref{condition_2}.  For $\rho > 0$, $\mcalS^{d_n-1}(\rho) = \{x \in \R^{d_n}: \|x\| = \rho\}$ denotes the hypersphere in $\R^{d_n}$ with center $0$ and radius $\rho$, and $\mcalS^{d_n-1}$ denotes the unit hypersphere $\mcalS^{d_n-1}(1)$.  We also use the notation $\|M\| := [\tr(MM')]^{1/2}$ for the Frobenius norm of any matrix $M$.

\begin{example}
\label{ex_Xn_Bingham}
{\rm 
Let $\{r_n, n \ge 1\}$ be a sequence of radii such that $r_n \to \sigma$ as $n \to \infty$, and suppose that the distribution of $X_n$ is concentrated on the hypersphere $\mcalS^{d_n-1}(r_n)$.  Since $\|X_n\|^2 = r_n^2$ then, trivially, $\|X_n\|^2 \cip \sigma^2$ and so \ref{condition_1} holds.  

Fix $\beta \in [0,1)$, and let $\{\Sigma_n, n \ge 1\}$ be a sequence of symmetric $d_n \times d_n$ matrices such that $\|\Sigma_n\| = O(d_n^{\beta/2})$ as $n \to \infty$.  Since $X_n \in \mcalS^{d_n-1}(r_n)$ then, by polar coordinates, $X_n = r_n \Theta_n$ where the random vector $\Theta_n \in \mcalS^{d_n-1}$.  Suppose also that $\Theta_n$ has a \textit{Bingham distribution with matrix parameter} $\Sigma_n$.  Relative to the surface measure $\dd\theta_n$ on $\mcalS^{d_n-1}$, normalized to have total surface area $1$, the probability density function of $\Theta_n$ is 
\begin{equation}
\label{eq_Bingham_pdf}
f(\theta;\Sigma_n) = [c(\Sigma_n)]^{-1} \exp(\theta' \Sigma_n \theta),
\end{equation}
$\theta \in \mcalS^{d_n-1}$.  
The normalizing constant $c(\Sigma_n)$ can be expressed in terms of a confluent hypergeometric function of matrix argument (cf., \cite{Bagyan_Richards_24}, \cite{Bingham}, or \citet[p.~288]{Muirhead}), however we will not need that result explicitly.  

It is simple to verify that for any $\tau \in \R$, $f(\theta;\Sigma_n - \tau I_{d_n}) \equiv f(\theta;\Sigma_n)$.  Therefore, with no loss of generality, we assume in \eqref{eq_Bingham_pdf} that $\tr(\Sigma_n) = 0$.  It is also evident that $\Theta_n \eqdist -\Theta_n$; therefore $\E(\Theta_n) = 0$ and hence $\E(X_n) = 0$.  Thus, with $\wtilX_n$ denoting an independent copy of $\wtilX_n$, we have $\E(X_n'\wtilX_n) = \E(X_n)'\E(\wtilX_n) = 0$.  

Next, observe that 
$$
\Var(X_n'\wtilX_n) = \E[(X_n'\wtilX_n)^2] 
= \E(X_n'\wtilX_n \cdot \wtilX_n'X_n) = \E \tr[(X_n X_n')(\wtilX_n \wtilX_n')].
$$
Interchanging the expectation and trace operations, and applying the independence of $X_n$ and $\wtilX_n$, we obtain the general identity,
\begin{equation}
\label{eq_Var_X_tilX}
\Var(X_n'\wtilX_n) = \tr [\E(X_nX_n') \E(\wtilX_n\wtilX_n')] 
= \tr\big([\Cov(X_n)]^2\big) = \big\|\Cov(X_n)\big\|^2,
\end{equation}
which is valid for any random vector $X_n$ and independent copy $\wtilX_n$ such that $\E(X_n) = 0$.  In the sequel, we will need to apply \eqref{eq_Var_X_tilX} repeatedly.  

Again resorting to polar coordinates $X_n = r_n \Theta_n$, the general identity \eqref{eq_Var_X_tilX} yields 
\begin{equation}
\label{eq_Var_X_Xtil_Bingham}
\Var(X_n'\wtilX_n) = \big\|\Cov(r_n \Theta_n)\big\|^2 = r_n^4 \, \big\|\Cov(\Theta_n)\big\|^2.
\end{equation}
Since $\|\Sigma_n\| = O(d_n^{\beta/2})$ as $n \to \infty$, where $\beta \in [0,1)$, then by \citet[Theorem 3.3, \textit{infra}]{Bagyan_Richards_24}, we obtain the expansion 
$$
\Cov(\Theta_n) = d_n^{-1} I_{d_n} + 2 d_n^{-1} (d_n+2)^{-1} \Sigma_n + O(d_n^{-(3-2\beta)/2}).
$$
On squaring both sides of this expansion, and recalling that $\tr(\Sigma_n) = 0$, we obtain 
\begin{equation}
\label{eq_E_Theta_Theta}
\big\|\Cov(\Theta_n)\big\|^2 = \tr[\Cov(\Theta_n)]^2 = d_n^{-1} + O\big(d_n^{-(3-2\beta)/2}\big);
\end{equation}
therefore $\big\|\Cov(\Theta_n)\big\|^2 \to 0$ as $n \to \infty$.  Since $r_n \to \sigma$ then it follows from \eqref{eq_Var_X_Xtil_Bingham} that $\Var(X_n'\wtilX_n) \to 0$, hence $X_n'\wtilX_n \cip 0$ as $n \to \infty$, so \ref{condition_2} holds. 

For the special case in which $\Sigma_n \equiv 0$, so that $X_n$ is uniformly distributed on $\mcalS^{d_n-1}(r_n)$, the above example was obtained by \citet[pp.~22--23]{Bagyan}.  
}\end{example}

\smallskip

In the next example, which was initiated by \citet[p.~23]{Bagyan}, we denote the $d_n$-dimensional ball centered at $0$ and radius $\rho$ by $\mcalB^{d_n}(\rho) = \{x \in \R^{d_n}: \|x\| \le \rho\}$ and the volume of the ball by $\textrm{Vol}(\mcalB^{d_n}(\rho))$.  

\begin{example}
\label{ex_Xn_ud_on_ball}
{\rm 
For a positive sequence $\{r_n, n \ge 1\}$ such that $r_n \to \sigma$ as $n \to \infty$, let $X_n$ be uniformly distributed on $\mcalB^{d_n}(r_n)$.  By polar coordinates, $X_n \eqdist R_n \Theta_n$ where $R_n \eqdist \|X_n\| \in [0,r_n]$, $\Theta_n$ is uniformly distributed on $\mcalS^{d_n-1}$, and $R_n$ and $\Theta_n$ are independent.  

Denote by $\dd\theta$ the normalized surface measure on $\mcalS^{d_n-1}$.  Applying polar coordinates on $\mcalB^{d_n}(r_n)$, viz., $x = s \theta$ where $0 \le s \le r_n$ and $\theta \in \mcalS^{d_n-1}$, together with the well-known formula for $\textrm{Vol}(\mcalB^{d_n}(1))$, we obtain 
\begin{align*}
\E \exp(\i t\|X_n\|^2) &= \frac{1}{\textrm{Vol}(\mcalB^{d_n}(1))} \int_{\mcalB^{d_n}(r_n)} \exp(\i t \|x\|^2) \dd x \nonumber \\
&= d_n r_n^{-d_n} \int_0^{r_n} s^{d_n - 1} \exp(\i t s^2) \dd s,
\end{align*}
Making the transformation $s \to r_n s^{1/d_n}$, we obtain 
$$
\E \exp(\i t\|X_n\|^2) 
= \int_0^1 \exp(\i t r_n^2 s^{2/d_n}) \dd s \to \int_0^1 \exp(\i t \sigma^2) \dd s 
= \exp(\i t \sigma^2)
$$
as $n \to \infty$.  Therefore $\|X_n\|^2 \cid \sigma^2$, hence $\|X_n\|^2 \cip \sigma^2$, so \ref{condition_1} holds.  

Let $\wtilX_n = \wtilde{R}_n \wtilde{\Theta}_n$ be an independent copy of $X_n$.  Since $\E(X_n) = 0$ then it follows that $\E(X_n'\wtilX_n) = 0$.  By the general identity \eqref{eq_Var_X_tilX}, 
\begin{equation}
\label{eq_Var_X_tilX_2}
\Var(X_n'\wtilX_n) = \big\|\Cov(X_n)\big\|^2 = \big\|\Cov(R_n \Theta_n)\big\|^2 = [\E(R_n^2)]^2 \, \big\|\Cov(\Theta_n)\big\|^2.
\end{equation}
We have $[\E(R_n^2)]^2 = [\E(\|X_n\|^2)]^2 \to \sigma^4$.  Also, by applying \eqref{eq_E_Theta_Theta} for the case in which $\Sigma_n = 0$, we obtain $\Var(X_n'\wtilX_n) \to 0$, so $X_n'\wtilX_n \cip 0$ as $n \to \infty$, hence \ref{condition_2} holds.  

This example can be extended further to the case in which $X_n$ has a dilated Bingham distribution, \textit{i.e.}, $X_n \eqdist R_n \Theta_n$ where $R_n$ is random; $R_n \cip \sigma$; $R_n$ and $\Theta_n$ are independent; $\Theta_n$ has a Bingham distribution with the density function \eqref{eq_Bingham_pdf}; and, as in Example \ref{ex_Xn_Bingham}, there exists $\beta \in [0,1)$ such that $\|\Sigma_n\| = O(d_n^{\beta/2})$ as $n \to \infty$.  In this setting, since $\|X_n\|^2 = R_n^2 \cip \sigma^2$ then \ref{condition_1} holds.  Also, proceeding as in \eqref{eq_Var_X_tilX_2}, we obtain $\E(X_n'\wtilX_n) = 0$ and $\Var(X_n'\wtilX_n) \to 0$.  Therefore $X_n'\wtilX_n \cip 0$, so \ref{condition_2} holds. 
}\end{example}

\begin{example}
\label{ex_Xn_ud_hyper_rect}
{\rm 
For $l_1,\ldots,l_n > 0$, set $\mcalC^{d_n}(l_n) = \{(x_1,\ldots,x_{d_n}) \in \R^{d_n}: |x_i| \le l_n/2, \, i=1,\ldots,d_n\}$, the $d_n$-dimensional hypercube centered at $0$ and with sides of length $l_n$.  

Let $\{L_n, n \ge 1\}$ be continuous random variables that satisfy $d_n L_n^2 \cip 12 \sigma^2$ as $n \to \infty$.  
Conditional on $L_n$, let $X_n = (X_{n;1},\ldots,X_{n;d_n})'$ be uniformly distributed on the hypercube $\mcalC^{d_n}(L_n)$; then $X_{n;1}|L_n,\ldots,X_{n;d_n}|L_n$ are mutually independent and identically uniformly distributed on the interval $[-L_n/2,L_n/2]$.  Therefore $\E(X_{n;1}|L_n) = 0$, $\E(X_{n;1}^2|L_n) = L_n^2/12$, and 
\begin{equation}
\label{eq_E_Xn_sq_Ln}
\E(\|X_n\|^2|L_n) = \E(X_{n;1}^2+\cdots+X_{n;d_n}^2|L_n) = d_n \E(X_{n;1}^2|L_n) = d_n L_n^2/12,
\end{equation}
and 
\begin{equation}
\label{eq_Var_Xn_sq_Ln}
\Var(\|X_n\|^2|L_n) = \sum_{j=1}^{d_n} \Var(X_{n;j}^2|L_n) = d_n \Var(X_{n;1}^2|L_n) = d_n L_n^4/180. 
\end{equation}
By \eqref{eq_E_Xn_sq_Ln} and the law of total expectation \citep[p.~333]{Ross}, 
$$
\E(\|X_n\|^2) = \E_{L_n} \E(\|X_n\|^2 | L_n) = \E_{L_n} (d_n L_n^2/12) \to \sigma^2.
$$
By \eqref{eq_E_Xn_sq_Ln}, \eqref{eq_Var_Xn_sq_Ln}, and the law of total variance \citep[p.~348]{Ross}, 
\begin{align}
\label{eq_total_var}
\Var(\|X_n\|^2) &= \E_{L_n} [\Var(\|X_n\|^2|L_n)] + \Var_{L_n} \big(\E(\|X_n\|^2|L_n)\big) \nonumber \\
&= (4/5) d_n^{-1} \E [(d_n L_n^2/12)^2] + \Var (d_n L_n^2/12).
\end{align}
Since $d_n L_n^2/12 \cip \sigma^2$ then $d_n^{-1} \E [(d_n L_n^2/12)^2] \to 0$ and $\Var [d_n L_n^2/12] \to 0$.  Therefore, by \eqref{eq_total_var}, $\Var(\|X_n\|^2) \to 0$ as $n \to \infty$, hence $\|X_n\|^2 \cip \sigma^2$ and \ref{condition_1} holds.  

Next, since $\E(X_n) = 0$ and $X_n$ and $\wtilX_n$ are independent then $\E(X_n'\wtilX_n) = 0$.  Also, it is simple to verify that $\Cov(X_n|L_n) = L_n^2 I_{d_n}/12$, hence $\tr\big([\Cov(X_n|L_n)]^2 = d_n (L_n^2/12)^2$.  Applying the general identity \eqref{eq_Var_X_tilX}, we obtain 
\begin{align*}
\Var(X_n'\wtilX_n) &= \E_{L_n} \tr\big([\Cov(X_n|L_n)]^2\big) \\ 
&= \E_{L_n} [d_n (L_n^2/12)^2] = d_n^{-1} \E_{L_n} [(d_n L_n^2/12)^2] \to 0.
\end{align*}
Therefore $X_n'\wtilX_n \cip 0$, so \ref{condition_2} holds.  

For the case in which the sequence $\{L_n, n \ge 1\}$ is deterministic, this example is due to \citet[p.~23]{Bagyan}.
}\end{example}

\begin{example}
\label{ex_Xn_gaussian}
{\rm 
Let $X_n \sim \mcalN_{d_n}(0,\Sigma_n)$ where $\Sigma_n$, the covariance matrix of $X_n$, is positive definite.  We suppose that $\tr(\Sigma_n) \to \sigma^2$ and $\tr(\Sigma_n^2) \to 0$ as $n \to \infty$.  

Denote by $\lambda_{n;1},\ldots,\lambda_{n;d_n}$ the eigenvalues of $\Sigma_n$, and let $H_n$ be a $d_n \times d_n$ orthogonal matrix such that $H_n \Sigma_n H_n'= \diag(\lambda_{n;1},\ldots,\lambda_{n;d_n})$.  Making the transformation $U_n = H_nX_n$ we find that $U_{n;1},\ldots,U_{n;d_n}$, the components of $U_n$, are mutually independent, with $U_{n;j} \eqdist \lambda_{n;j}^{1/2} Z_{n;j}$ with $Z_{n;1},\ldots,Z_{n;d_n}$ being mutually independent $\mcalN_1(0,1)$ random variables.  Therefore $\|X_n\|^2 = \|U_n\|^2 \eqdist \sum_{j=1}^{d_n} \lambda_{n;j} Z_{n;j}^2$, and it follows easily that 
$\E(\|X_n\|^2) 
= \tr(\Sigma_n)$, hence $\E(\|X_n\|^2) \to \sigma^2$ as $n \to \infty$.  Further, 
$\Var(\|X_n\|^2) 
= 2 \tr(\Sigma_n^2)$, 
so $\Var(\|X_n\|^2) \to 0$.  Therefore $\|X_n\|^2 \cip \sigma^2$ as $n \to \infty$, so \ref{condition_1} holds.  

For $\wtilX_n$, an independent copy of $X_n$, we have $\E(X_n'\wtilX_n) = 0$.  Applying the general identity \eqref{eq_Var_X_tilX} we find that \ref{condition_2} holds since, as $n \to \infty$, 
$$
\Var(X_n'\wtilX_n) = \tr\big([\Cov(X_n)]^2\big) = \tr(\Sigma_n^2)  \to 0.
$$

We now present two examples of $\Sigma_n$ such that $\tr(\Sigma_n) \to \sigma^2$ and $\tr(\Sigma_n^2) \to 0$.  For the first such example, suppose that 
\begin{equation}
\label{eq_Xn_gaussian_2}
\lambda_{n;j} = \sigma^2 (\log d_n)^{-1} j^{-1},
\end{equation}
$j=1,\ldots,d_n$.  Let $\gamma = 0.57721\ldots$ denote Euler's constant, which arises in the asymptotic formula \citep[\href{https://dlmf.nist.gov/2.10}{\S 2.10}]{Olver_Wong}, 
\begin{equation}
\label{eq_euler-maclaurin}
\sum_{j=1}^{d_n} j^{-1} = \gamma + \log d_n + O(d_n^{-1})
\end{equation}
as $n \to \infty$.  By \eqref{eq_euler-maclaurin}, we obtain 
$\tr(\Sigma_n) = \sigma^2 (\log d_n)^{-1} [\gamma + \log d_n + O(d_n^{-1})] \to \sigma^2$.  
Since $\sum_{j=1}^{d_n} j^{-2} < \pi^2/6$, then we also have 
$\tr(\Sigma_n^2) 
< \sigma^4 (\log d_n)^{-2} \pi^2/6 \to 0$.  

For the second example, let $r > -1/2$ and define
\begin{equation}
\label{eq_Xn_gaussian_3}
\lambda_{n;j} = (r+1) \sigma^2 d_n^{-(r+1)} j^r,
\end{equation}
$j=1,\ldots,d_n$.  The distinction between \eqref{eq_Xn_gaussian_2} and \eqref{eq_Xn_gaussian_3} is that, for fixed $n$, \eqref{eq_Xn_gaussian_2} is decreasing in $j$ whereas \eqref{eq_Xn_gaussian_3} is increasing in $j$.  Applying Euler-Maclaurin summation, 
we have 
\begin{equation}
\label{eq_Euler_Maclaurin}
\sum_{j=1}^{d_n} j^r 
= (r+1)^{-1} d_n^{r+1} [1 + O(d_n^{-1})].
\end{equation}
Letting $n \to \infty$, it follows from \eqref{eq_Xn_gaussian_3} and \eqref{eq_Euler_Maclaurin} that  
$\tr(\Sigma_n) 
= [1 + O(d_n^{-1})] \sigma^2 \to \sigma^2$ 
and 
$\tr(\Sigma_n^2) 
= (r+1)^2 (2r+1)^{-1} \sigma^4 d_n^{-1} [1 + O(d_n^{-1})] \to 0$.  
}\end{example}

\smallskip

Next, we provide an example for which \ref{condition_1} does not hold whereas \ref{condition_2} holds.

\begin{example}
\label{ex_Xn_t_distn}
{\rm 
For $\nu > 4$, let $X_n$ have a centered multivariate $t$-distribution with index parameter $\nu$ and positive definite matrix parameter $\Sigma_n$ \citep[p.~48]{Muirhead}.  There holds the stochastic representation $X_n \eqdist \nu^{1/2} Q_\nu^{-1/2} Z_n$ where $Q_\nu \sim \chi^2_{\nu}$, a chi-squared distribution with $\nu$ degrees-of-freedom, $Z_n \sim \mcalN_{d_n}(0,\Sigma_n)$, and $Q_\nu$ and $Z_n$ are independent.  We also assume that $\tr(\Sigma_n) \to (\nu-2)\sigma^2/\nu$ and $\tr(\Sigma_n^2) \to 0$ as $n \to \infty$.  

It is straightforward to verify that $\E(X_n) = 0$ and that $\Cov(X_n) = \E(X_nX_n') = \nu\Sigma_n/(\nu-2)$.  Also, $\E(X_n'\wtilX_n) = 0$ and, by the general identity \eqref{eq_Var_X_tilX}, 
$$
\Var(X_n'\wtilX_n) = \tr \big([\Cov(X_n)]^2\big) = \frac{\nu^2}{(\nu-2)^2} \tr(\Sigma_n^2).
$$
Therefore $\Var(X_n'\wtilX_n) \to 0$ as $n \to \infty$, so $X_n'\wtilX_n \cip 0$ and \ref{condition_2} holds.   

In considering \ref{condition_1}, we begin by noting that 
$$
\E(\|X_n\|^2) 
= \tr \big(\Cov(X_n)\big) = \frac{\nu}{\nu-2}\tr(\Sigma_n) \to \sigma^2,
$$
as $n \to \infty$.  Applying the law of total variance \citep[p.~348]{Ross}, and the independence of $Q_\nu$ and $Z_n$, we obtain
\begin{align*}
\Var(\|X_n\|^2) 
&= \nu\Big(\E(Q_\nu^{-2}) \Var(\|Z_n\|^2) + \Var(Q_\nu^{-1}) [\E(\|Z_n\|^2)]^2\Big) \nonumber \\
&\ge \nu\Var(Q_\nu^{-1}) [\E(\|Z_n\|^2)]^2 \nonumber \\
&= \frac{\nu}{(\nu-4)(\nu-2)^2} [\tr(\Sigma_n)]^2,
\end{align*}
and therefore 
$$
\lim_{n \to \infty} \Var(\|X_n\|^2) \ge \frac{\nu}{(\nu-4)(\nu-2)^2} \lim_{n \to \infty} [\tr(\Sigma_n)]^2 
= \frac{\sigma^4}{\nu(\nu-4)} > 0.
$$
Since 
$$
\lim_{n \to \infty} \E(\|X_n\|^4) = \lim_{n \to \infty} \Var(\|X_n\|^2) + \lim_{n \to \infty} [\E(\|X_n\|^2)]^2 \ge \frac{\sigma^4}{\nu(\nu-4)} + \sigma^4 > \sigma^4,
$$
then it follows that $\|X_n\|^4 \not\cip \sigma^4$.  Therefore $\|X_n\|^2 \not\cip \sigma^2$, so \ref{condition_1} does not hold.  

To complete this example, we note that the Laplace distributions also satisfy \ref{condition_2} but not \ref{condition_1}.  For those distributions, $X_n \eqdist Q_\nu^{1/2} Z_n$ where $Q_\nu \sim \chi^2_{\nu}$, $Z_n \sim \mcalN_{d_n}(0,\Sigma_n)$, and $Q_\nu$ and $Z_n$ are mutually independent.  
}\end{example}

\section{Properties of the probability density function of \texorpdfstring{$\boldsymbol{Y_n|\Xi_n}$}{YngivenXin}}
\label{sec_Lp_conv_cond_pdf}

\subsection{Preliminary remarks on the vectors \texorpdfstring{$\boldsymbol{X_n}$}{Xn} and \texorpdfstring{$\boldsymbol{\Xi_n}$}{Xin}}
\label{subsec_prelim_remarks}

Let $\{X_n \in \R^{d_n}, n \ge 1\}$ be a sequence of continuous random vectors, each satisfying \ref{condition_1} and \ref{condition_2}.  We assume that the random modulators $\{\Xi_n \in \R^{d_n}, n \ge 1\}$ are continuous and mutually independent of $\{X_n, n \ge 1\}$.  Also denote by $f_{X_n}$ and $f_{\Xi_n}$ the marginal density functions of $X_n$ and $\Xi_n$, respectively, each density assumed to being supported on an open subset of $\R^{d_n}$.  

Let $Y_n = \Xi_n'X_n$; then we obtain the joint density function of $(Y_n,\Xi_n)$ by making the usual transformation from $(X_n,\Xi_n)$ to $(Y_n,X_{n;2},\ldots,X_{n;d_n},\Xi_n)$, where $X_{n;j}$ is the $j$th component of $X_n$, $j=2,\ldots,d_n$.  
Since $\Xi_n$ is continuous then the Jacobian of the transformation exists and is non-zero, almost everywhere.  Therefore $f_{Y_n,\Xi_n}$, the joint density function of $(Y_n,\Xi_n)$, exists almost everywhere and is obtained by integrating over the support of $X_{n;2},\ldots,X_{n;d_n}$.  Consequently $f_{Y_n|\Xi_n}$, the conditional density function of $Y_n|\Xi_n$, also exists almost everywhere and 
$$
f_{Y_n|\Xi_n}(y) = f_{Y_n,\Xi_n}(y,\xi)/f_{\Xi_n}(\xi)
$$
for all $y \in \R$ and all $\xi \in \R^{d_n}$ such that $f_{\Xi_n}(\xi) \neq 0$.  

We will also encounter the conditional characteristic function of $Y_n|\Xi_n$, viz., 
$$
\varphi_{Y_n|\Xi_n}(t) = \E_{Y_n|\Xi_n} \exp(\i tY_n) = \E_{X_n|\Xi_n} \exp(\i t\Xi_n'X_n),
\quad t \in \R.
$$
The following result provides in terms of $\varphi_{X_n}$, the characteristic function of $X_n$, a condition under which $\varphi_{Y_n|\Xi_n}$ is integrable for almost all values of $\Xi_n$.  

\begin{lemma}
\label{lem_integrable_cf_Yn_Xin}
A necessary and sufficient condition that $\varphi_{Y_n|\Xi_n} \in L^1(\R)$ for almost all values of $\Xi_n$ is that, for almost all $\theta \in \mcalS^{d_n-1}$, 
\begin{equation}
\label{eq_integrable_cf_Yn_Xin}
\int_{-\infty}^\infty \big|\varphi_{X_n}(t\theta)\big| \dd t < \infty.
\end{equation}
\end{lemma}

We now assume that the distribution of $\Xi_n$ is orthogonally invariant, \textit{i.e.}, $\Xi_n \eqdist H \Xi_n$ for all $d_n \times d_n$ orthogonal matrices $H$ \citep[p.~34]{Muirhead}.  It is well known that this orthogonal invariance is equivalent to the property that $\Xi_n$ has a spherically symmetric characteristic function, \textit{i.e.}, $\E \exp(\i u'\Xi_n)$, $u \in \R^{d_n}$, depends on $\|u\|$ only.  We assume that there exists a function $\psi: [0,\infty) \to \R$ such that, for all $n = 1,2,3,\ldots$, 
\begin{equation}
\label{eq_Xi_n_cf}
\E_{\Xi_n} \exp(\i  u'\Xi_n) = \psi(\|u\|^2), \quad u \in \R^{d_n}.
\end{equation}
By a famous theorem of \citet{Schoenberg} (see also \citet{Eaton}, \citet{Ressel}, \citet{Steerneman}), there exists a distribution function $G$ on $[0,\infty)$ such that 
\begin{equation}
\label{eq_Schoenberg}
\psi(t^2) = \int_0^\infty \exp(- t^2 v^2/2) \dd G(v), 
\quad t \in \R,
\end{equation}
equivalently, $\psi(t^2)$ is a scale mixture of one-dimensional Gaussian characteristic functions.  

We require that the characteristic function \eqref{eq_Xi_n_cf} be integrable on $\R^{d_n}$.  By integration using polar coordinates, we have 
\begin{equation}
\label{eq_Schoenberg_15}
\int_{\R^{d_n}} \psi(\|u\|^2) \dd u \propto \int_0^\infty t^{d_n-1} \psi(t^2) \dd t,
\end{equation}
so the integrability requirement on \eqref{eq_Xi_n_cf} is equivalent to the function $t \mapsto t^{d_n-1} \psi(t^2)$, $t > 0$, being integrable.  On substituting in the right-hand side of \eqref{eq_Schoenberg_15} the expression for $\psi(t^2)$ from \eqref{eq_Schoenberg}, and then interchanging the order of integration in $t$ and $v$, it follows by Fubini's theorem that 
\begin{align*}
\int_0^\infty t^{d_n-1} \psi(t^2) \dd t 
&= \int_0^\infty \int_0^\infty t^{d_n-1} \exp(- t^2 v^2/2) \dd t \dd G(v) \\
&= 2^{(d_n-2)/2} \Gamma(d_n/2) \int_0^\infty v^{-d_n} \dd G(v).
\end{align*}
Therefore the integrability assumption on $\psi(\cdot)$ is equivalent to 
\begin{equation}
\label{eq_G_dn_integrable}
\int_0^\infty v^{-d_n} \dd G(v) < \infty.
\end{equation}
By \eqref{eq_Xi_n_cf} and \eqref{eq_Schoenberg}, 
\begin{equation}
\label{eq_Schoenberg_3}
\E_{\Xi_n} \exp(\i  u'\Xi_n) = \int_0^\infty \exp(- \|u\|^2 v^2/2) \dd G(v),
\quad u \in \R^{d_n}.
\end{equation}
Also applying to \eqref{eq_Schoenberg_3} the multidimensional inverse Fourier transform, it follows that the density function of $\Xi_n$ exists and is given by 
$$
f_{\Xi_n}(\xi) = \int_0^\infty (2\pi)^{-d_n/2} v^{-d_n} \exp(-\|\xi\|^2/2v^2) \dd G(v), 
\quad \xi \in \R^{d_n}.
$$

\subsection{The matrices \texorpdfstring{$\boldsymbol{A_{n,k}}$}{An}}
\label{subsec_A_nk_matrix}

For here on, we denote by $k$ a fixed integer.  Let $\wtilX_{n,1},\ldots,\wtilX_{n,k}$ be mutually independent copies of $X_n$; in particular, $\wtilX_{n,1},\ldots,\wtilX_{n,k}$ satisfy \ref{condition_1} and \ref{condition_2}.  
Define the $d_n \times k$ matrix $\wtilde{\mcalX}_n = (\wtilX_{n,1} \ \cdots \ \wtilX_{n,k})$ and the $k \times k$ positive semidefinite matrix 
\begin{equation}
\label{eq_A_nk_matrix}
A_{n,k} = \wtilde{\mcalX}_n' \wtilde{\mcalX}_n = \big(\wtilX_{n,j}'\wtilX_{n,r}\big)_{j,r=1}^k \, .
\end{equation}

We assume henceforth the following condition on the distribution of $X_n$: 
\begin{enumerate}[label=(C.\arabic*)]\addtocounter{enumi}{2}
\item \label{condition_3}
There exist a positive integer $n_0$ such that $\E[(\det A_{n_0,k})^{-1/2}] < \infty$.
\end{enumerate}

As a consequence of \ref{condition_3}, there holds the following properties of $A_{n,k}$:

\begin{lemma}
\label{lem_An_props}
Suppose that \ref{condition_3} holds.  Then for all $n \ge n_0$, 
\begin{enumerate}[label=(\roman*),nosep]
\item
$d_n \ge k$.
\item
$A_{n,k}$ is positive definite, almost surely.
\item
$\E[(\det A_{n+1,k})^{-1/2}] \le \E[(\det A_{n,k})^{-1/2}]$.
\item
$\E[(\det A_{n,j})^{-1/2}] < \infty$ for all $j=1,\ldots,k$.  
\end{enumerate}
\end{lemma}

\subsection{Convergence properties of the probability density function} 
\label{subsec_cim_pdf}

For $f \in L^1(\R)$, we introduce the notation 
\begin{equation}
\label{eq_Fourier_notat}
{\mcalF}_{\goesto{y}{t}} f(y) \equiv (\mcalF f)(t) = \int_{-\infty}^\infty \exp(\i ty) f(y) \dd y, 
\quad t \in \R, 
\end{equation}
for the Fourier transform of $f$.  For a Fourier transform $\widehat{f} \in L^1(\R)$, we often write 
\begin{equation}
\label{eq_Fourier_inv_notat}
\mcalF_{\goesto{t}{y}}^{-1} \widehat{f}(t) \equiv (\mcalF^{-1} \widehat{f}\,)(y) = (2\pi)^{-1} \int_{-\infty}^\infty \exp(-\i yt) \widehat{f}(t) \dd t, \quad y \in \R,
\end{equation}
for the inverse Fourier transform of $\widehat{f}$.  The notations ${\mcalF}_{\goesto{y}{t}}$ and $\mcalF_{\goesto{t}{y}}^{-1}$ will be used repeatedly to monitor the arguments of numerous simultaneous Fourier and inverse Fourier transforms, and we also use similar notation in fewer instances for the multidimensional Fourier and inverse Fourier transforms.  

We now state the main result of this section.  In this result and hereafter, $G$ denotes the distribution function defined by \eqref{eq_Schoenberg} and $V$ denotes the corresponding random variable; we use the notation 
$$
f_{\mcalN_k(0,\Sigma)}(w) = (2\pi)^{-k/2} (\det \Sigma)^{-1/2} \exp(- \tfrac12 w'\Sigma^{-1}w), 
\quad w \in \R^k,
$$
for the probability density function of the $k$-dimensional normal distribution with mean $0$ and covariance matrix $\Sigma$; and $\bfone_k$ denotes the vector $(1,\ldots,1)' \in \R^k$.

\begin{theorem}
\label{thm_Lp_conv_cond_pdf_spher}
Suppose that the random vectors $\{X_n \in \R^{d_n}, n \ge 1\}$ satisfy \ref{condition_1}, \ref{condition_2}, \ref{condition_3}, and \eqref{eq_integrable_cf_Yn_Xin}.  Let $\{\Xi_n \in \R^{d_n}, n \ge 1\}$ be spherically symmetric modulating vectors that satisfy \eqref{eq_Xi_n_cf} and \eqref{eq_G_dn_integrable} and are independent of $\{X_n, n \ge 1\}$, and let $Y_n = \Xi_n'X_n$, $n \ge 1$.  Then for all $y \in \R$ and all $j=1,\ldots,k$, 
\begin{equation}
\label{eq_conv_cond_pdf_spher}
\lim_{n \to \infty} \E_{\Xi_n} \big[f_{Y_n|\Xi_n}(y)\big]^j = \E_V \big[f_{\mcalN_1(0,\sigma^2 V^2)}(y)\big]^j.
\end{equation}
\end{theorem}

\medskip

For the case in which $\Xi_n \sim \mcalN_{d_n}(0,I_{d_n})$, it follows from \eqref{eq_Schoenberg_3} that $G$ is concentrated at $v=1$, hence \eqref{eq_G_dn_integrable} holds trivially.  Then we obtain the following $p$th-mean pointwise convergence property of $f_{Y_n|\Xi_n}$.  

\smallskip

\begin{corollary}
\label{cor_Lp_conv_cond_pdf}
Let $\{X_n \in \R^{d_n}, n \ge 1\}$ be continuous random vectors that satisfy \ref{condition_1}, \ref{condition_2}, \ref{condition_3}, and \eqref{eq_integrable_cf_Yn_Xin}, and let $\:\Xi_n \sim \mcalN_{d_n}(0,I_{d_n})$.  Then for all $\:y \in \R$ and $0 \le p < 2 \lfloor k/2 \rfloor$, 
\begin{equation}
\label{eq_Lp_conv_cond_pdf}
\lim_{n \to \infty} \E_{\Xi_n} \big|f_{Y_n|\Xi_n}(y) - f_{\mcalN_1(0,\sigma^2)}(y)\big|^p = 0.
\end{equation}
\end{corollary}

The following result quantifies explicitly a rate of convergence in \eqref{eq_conv_cond_pdf_spher} in terms of the regularity assumptions \ref{condition_1} and \ref{condition_2}, and therefore strengthens Theorem \ref{thm_Lp_conv_cond_pdf_spher}.  Moreover \eqref{eq_quant_pdf}, together with the calculations in Section \ref{sec_examples}, provides a rate of convergence for each example in that section.

\begin{theorem}
\label{thm_quant_pdf}
Suppose that $X_n$ and $\Xi_n$ satisfy the assumptions of Theorem \ref{thm_Lp_conv_cond_pdf_spher}, and let $1 \le j \le k$.  Then there exists $n_j \in \N$ such that, for all $n \ge n_j$, 
\begin{equation}
\label{eq_quant_pdf}
\sup_{y \in \R} \big|\E_{\Xi_n} \big[f_{Y_n|\Xi_n}(y)\big]^j - \E_V \big[f_{\mcalN_1(0,\sigma^2 V^2)}(y)\big]^j\big| \le c_j \, \big[\E \big\|A_{n,j} - \sigma^2 I_j\big\|^2\big]^{1/2},
\end{equation}
where 
\begin{equation}
\label{eq_quant_c_k}
c_j = 2^{-(j+1)/2} \pi^{-j/2} j^{3/2} \sigma^{-(j+2)} \E(V^{-j}).
\end{equation}
Further, 
\begin{equation}
\label{eq_quant_rate}
\E \big\|A_{n,j} - \sigma^2 I_j\big\|^2 = j \, \E (\|X_n\|^2 - \sigma^2)^2 + j(j-1) [\E (X_n'\wtilX_n)]^2.
\end{equation}
\end{theorem}

\subsection{Proofs}
\label{subsec_pdf_conv_proofs}

\noindent
\textit{Proof of Lemma \ref{lem_integrable_cf_Yn_Xin}}:  
Since $X_n$ and $\Xi_n$ are independent then, for all $t \in \R$ and $\xi \in \R^{d_n}$, 
$$
\varphi_{Y_n|\{\Xi_n=\xi\}}(t) = \E_{X_n} \exp(\i t\xi'X_n) = \varphi_{X_n}(t\xi).
$$
Therefore for $\xi \neq 0$, 
$$
\big\|\varphi_{Y_n|\{\Xi_n=\xi\}}\big\|_{L^1(\R)} := \int_{-\infty}^\infty |\varphi_{Y_n|\{\Xi_n=\xi\}}(t)| \dd t = \int_{-\infty}^\infty |\varphi_{X_n}(t\xi)| \dd t.
$$
Making the change-of-variable $t \mapsto t/\|\xi\|$, which is permissible since $\xi \neq 0$, we obtain 
$$
\big\|\varphi_{Y_n|\{\Xi_n=\xi\}}\big\|_{L^1(\R)} 
= \frac{1}{\|\xi\|} \int_{-\infty}^\infty |\varphi_{X_n}(t\xi/\|\xi\|)| \dd t 
= \frac{1}{\|\xi\|} \int_{-\infty}^\infty |\varphi_{X_n}(t\theta)| \dd t,
$$
where $\theta = \xi/\|\xi\| \in \mcalS^{d_n - 1}$.  Since the mapping $\xi \to \theta = \xi/\|\xi\|$ from $\R^{d_n}\setminus\{0\}$ to $\mcalS^{d_n-1}$ is surjective then it follows that $\varphi_{Y_n|\Xi_n} \in L^1(\R)$ if and only if \eqref{eq_integrable_cf_Yn_Xin} holds.  
$\qed$

\medskip

\noindent
\textit{Proof of Lemma \ref{lem_An_props}}:  
(i) Let $\wtilX_{n,j;1},\ldots,\wtilX_{n,j;d_n}$ be the components of $\wtilX_{n,j}$, $j=1,\ldots,k$, then by \eqref{eq_A_nk_matrix}, 
\begin{align}
\label{eq_An_decomp}
A_{n,k} = \bigg(\sum_{m=1}^{d_n} \wtilX_{n,j;m} \wtilX_{n,r;m}\bigg)_{j,r=1}^k 
&\equiv \sum_{m=1}^{d_n} 
\begin{pmatrix} \wtilX_{n,1;m} \\ \vdots \\ \wtilX_{n,k;m} \end{pmatrix} 
(\wtilX_{n,1;m},\ldots,\wtilX_{n,k;m}),
\end{align}
which is a sum of $d_n$ positive semidefinite matrices.  By \ref{condition_3}, $A_{n_0,k}$ is nonsingular, almost surely, so it follows by \eqref{eq_An_decomp} that $d_{n_0} \ge k$.  Since $\{d_n, n \ge 1\}$ is increasing then we have $d_n \ge k$ for all $n \ge n_0$.  

(ii)  By \eqref{eq_An_decomp}, $A_{n,k}$ is positive semidefinite, so $\det(A_{n,k}) \ge 0$.  Therefore to prove that $A_{n,k}$ is positive definite (almost surely), it suffices to show that $\det(A_{n,k}) > 0$, almost surely.  

It is evident that $\wtilde{\mcalX}_n$ has a probability density function on the underlying Euclidean space $\R^{d_n k}$.  Therefore, by a result of \citet[p.~344]{Malley}, the probability distribution of $\wtilde{\mcalX}_n$ assigns zero probability to the zeros of any non-trivial polynomial in the components of $\wtilde{\mcalX}_n$.  Since $\det(A_{n,k})$ is a non-trivial polynomial in the components of $\wtilde{\mcalX}_n$ then, by Malley's theorem, $\P\big(\det(A_{n,k}) = 0\big) = 0$.  Therefore $\det(A_{n,k}) > 0$, almost surely.  

(iii) Since $d_{n+1} \ge d_n$ then, by \eqref{eq_An_decomp}, $A_{n+1,k} - A_{n,k}$ is positive semidefinite.  By \citet[p.~495, Corollary 7.7.4(e)]{Horn} we obtain $\det(A_{n+1,k}) \ge \det(A_{n,k})$,  hence $\E[(\det A_{n+1,k})^{-1/2}] \le \E[(\det A_{n,k})^{-1/2}]$.  

(iv) Since $A_{n,k}$ is positive semidefinite then, by Hadamard's inequality \citep[p.~505]{Horn}, 
$$
\det A_{n,k} \le \prod_{j=1}^k \wtilX_{n,j}'\wtilX_{n,j} = \prod_{j=1}^k \|\wtilX_{n,j}\|^2.
$$
Since $\wtilX_{n,1},\ldots,\wtilX_{n,k}$ are mutually independent copies of $X_n$ then it follows that 
$$
\E \big[(\det A_{n,k})^{-1/2}\big] \ge \E \prod_{j=1}^k \|\wtilX_{n,j}\|^{-1} = \big(\E \|X_n\|^{-1}\big)^k.
$$
As shown before, $\E \big[(\det A_{n,k})^{-1/2}\big] < \infty$, so it follows that $\E(\|X_n\|^{-1}) < \infty$.  Noting that $\|X_n\|^2 = A_{n,1}$ then we have also shown that $\E \big[(\det A_{n,1})^{-1/2}\big] < \infty$.

Define the vector $D_{n,k-1} = (\wtilX_{n,1}'\wtilX_{n,k},\ldots,\wtilX_{n,k-1}'\wtilX_{n,k})'$; then we can write $A_{n,k}$ in partitioned form, 
$$
A_{n,k} = 
\begin{pmatrix}
A_{n,k-1}  & D_{n,k-1} \\
D_{n,k-1}' & \wtilX_{n,k}'\wtilX_{n,k}
\end{pmatrix}.
$$
Since $A_{n,k}$ is positive semidefinite then, by the Hadamard-Fischer inequality \citep[p.~506]{Horn}, 
$$
\det(A_{n,k}) \le \det(A_{n,k-1}) \cdot (\wtilX_{n,k}'\wtilX_{n,k}) = \|\wtilX_{n,k}\|^2 \, \det(A_{n,k-1}),
$$
equivalently, 
$$
\|\wtilX_{n,k}\|^{-1} (\det A_{n,k-1})^{-1/2} \le (\det A_{n,k})^{-1/2}.
$$
As $\wtilX_{n,k}$ is independent of $A_{n,k-1}$ and since $\E[(\det A_{n,k})^{-1/2}] < \infty$ then, by taking expectations, we obtain 
$$
\E(\|\wtilX_{n,k}\|^{-1}) \cdot \E[(\det A_{n,k-1})^{-1/2}] \le \E[(\det A_{n,k})^{-1/2}] < \infty.
$$
Therefore $\E[(\det A_{n,k-1})^{-1/2}] < \infty$.  By repeating this argument, we deduce finally that $\E[(\det A_{n,j})^{-1/2}] < \infty$ for all $j=k-1,k-2,\ldots,2$.  
$\qed$

\medskip

\noindent
\textit{Proof of Theorem \ref{thm_Lp_conv_cond_pdf_spher}}:  
Consider the case in which $j=k$.  Applying the Fourier transform with the notation \eqref{eq_Fourier_notat}, we have 
$$
\varphi_{Y_n|\Xi_n}(t) = \E_{Y_n|\Xi_n} \exp(\i tY_n) 
= \mcalF_{\goesto{y}{t}} f_{Y_n|\Xi_n}(y),
$$
$t \in \R$.  By \eqref{eq_integrable_cf_Yn_Xin}, $\varphi_{Y_n|\Xi_n}$ is integrable, so by applying \eqref{eq_Fourier_inv_notat} to invert the Fourier transform of $f_{Y_n|\Xi_n}$, it follows that, for all $y \in \R$,  
\begin{equation}
\label{eq_pdf_Yn_Xin}
f_{Y_n|\Xi_n}(y) = \mcalF_{\goesto{t}{y}}^{-1} \varphi_{Y_n|\Xi_n}(t) = \mcalF^{-1}_{\goesto{t}{y}} \E_{X_n|\Xi_n} \exp(\i t\Xi_n'X_n).
\end{equation}

Since $\wtilX_{n,1},\ldots,\wtilX_{n,k}$ are mutually independent copies of $X_n$ then, by \eqref{eq_pdf_Yn_Xin}, 
$$
\big[f_{Y_n|\Xi_n}(y)\big]^k 
= \prod_{j=1}^k \mcalF_{\goesto{t_j}{y}}^{-1} \E_{\wtilX_{n,j}} \exp(\i t_j \Xi_n'\wtilX_{n,j}).
$$
After formally interchanging expectations and inverse Fourier transforms, we obtain 
\begin{align}
\label{eq_converg_expect_cond_pdf_2}
\E_{\Xi_n} \big[f_{Y_n|\Xi_n}(y)\big]^k &= \E_{\Xi_n} \prod_{j=1}^k \mcalF_{\goesto{t_j}{y}}^{-1} \E_{\wtilX_{n,j}} \exp(\i  t_j \Xi_n'\wtilX_{n,j}) \nonumber \\
&= \Big(\prod_{j=1}^k \mcalF_{\goesto{t_j}{y}}^{-1} \, \E_{\wtilX_{n,j}}\Big) \E_{\Xi_n} \exp\Big(\i \Xi_n' \sum_{j=1}^k t_j \wtilX_{n,j}\Big).
\end{align}

Let $w = (t_1,\ldots,t_k)' \in \R^k$.  Since $\Xi_n$ is spherically symmetric with characteristic function \eqref{eq_Xi_n_cf} then, conditional on $\wtilde{\mcalX}_n$, 
\begin{equation}
\label{eq_Ank_psi}
\E_{\Xi_n|\wtilde{\mcalX}_n} \exp\Big(\i \Xi_n' \sum_{j=1}^k t_j \wtilX_{n,j}\Big) = \psi\Big(\Big\|\sum_{j=1}^k t_j \wtilX_{n,j}\Big\|^2\Big) \equiv \psi(w' A_{n,k} w).
\end{equation}
Substituting this result in \eqref{eq_converg_expect_cond_pdf_2} and again formally interchanging Fourier transforms and expectations, we obtain  
\begin{align*}
\E_{\Xi_n} \big[f_{Y_n|\Xi_n}(y)\big]^k &= \Big(\prod_{j=1}^k \mcalF_{\goesto{t_j}{y}}^{-1} \E_{\wtilX_{n,j}}\Big) \psi(w' A_{n,k} w) \nonumber \\
&= \E_{\wtilX_{n,1}} \cdots \E_{\wtilX_{n,k}} \mcalF_{\goesto{t_1}{y}}^{-1} \cdots \mcalF_{\goesto{t_k}{y}}^{-1} \psi(w' A_{n,k} w).
\end{align*}
Since $w = (t_1,\ldots,t_k)'$ then a moment of reflection reveals that 
$$
\mcalF_{\goesto{t_1}{y}}^{-1} \cdots \mcalF_{\goesto{t_k}{y}}^{-1} \equiv \mcalF_{\goesto{w}{y \bfone_k}}^{-1},
$$
the $k$-dimensional inverse Fourier transform, evaluated at $y\bfone_k$, of a function of $w$.  Therefore 
\begin{equation}
\label{eq_Fourier_inv}
\E_{\Xi_n} \big[f_{Y_n|\Xi_n}(y)\big]^k = \E_{\wtilde{\mcalX}_n} \mcalF_{\goesto{w}{y \bfone_k}}^{-1} \psi(w' A_{n,k} w).
\end{equation}

Recall that, for $\widehat{f} \in L^1(\R^k)$, the inverse Fourier transform is 
\begin{equation}
\label{eq_FT_inv_integral}
(\mcalF^{-1} \widehat{f}\:)(u) \equiv 
\mcalF_{\goesto{w}{u}}^{-1} \widehat{f}(w) = (2\pi)^{-k} \int_{\R^k} \exp(-\i u'w) \, \widehat{f}(w) \dd w, \quad u \in \R^k,
\end{equation}
so 
\begin{equation}
\label{eq_FT_inv_integral_2}
\mcalF_{\goesto{w}{u}}^{-1} \psi(w' A_{n,k} w) = (2\pi)^{-k} \int_{\R^k} \exp(-\i u'w) \, \psi(w' A_{n,k} w) \dd w, \quad u \in \R^k.
\end{equation}
By \eqref{eq_G_dn_integrable}, $t^{d_n-1} \psi(t^2) \in L^1(0,\infty)$; since $d_n \ge k$ then, by Jensen's inequality, $t^{k-1} \psi(t^2) \in L^1(0,\infty)$.  Arguing as in \eqref{eq_Schoenberg_15}-\eqref{eq_G_dn_integrable}, we also obtain $\psi(\|w\|^2) \in L^1(\R^k)$, and therefore $\psi(w'Aw) \in L^1(\R^k)$ for positive definite $A$; hence \eqref{eq_FT_inv_integral_2} is finite for all $w$.  On substituting for $\psi(\cdot)$ from \eqref{eq_Schoenberg} and formally interchanging integrals, we obtain 
\begin{equation}
\label{eq_Fourier_inv_15}
\mcalF_{\goesto{w}{u}}^{-1} \psi(w' A_{n,k} w) 
= (2\pi)^{-k} \int_0^\infty \int_{\R^k} \exp(-\i u'w - \tfrac12 v^2 w'A_{n,k}w) \dd w \dd G(v),
\end{equation}
$u \in \R^k$.  By Lemma \ref{lem_An_props}(ii), $A_{n,k}$ is nonsingular, almost surely, for $d_n \ge k$; then by applying to \eqref{eq_Fourier_inv_15} the multivariate Gaussian integral, \textit{viz.}, 
$$
\int_{\R^k} \exp(-\i u'w - \tfrac12 v^2 w'A_{n,k}w) \dd w 
= (2\pi)^{k/2} v^{-k} (\det A_{n,k})^{-1/2} \exp(-\tfrac12 v^{-2} u'A_{n,k}^{-1} u),
$$
and simplifying the resulting expression, we obtain for $d_n \ge k$, 
\begin{equation}
\label{eq_Fourier_inv_18}
\mcalF_{\goesto{w}{u}}^{-1} \psi(w' A_{n,k} w) = (2\pi)^{-k/2} (\det A_{n,k})^{-1/2} \int_0^\infty v^{-k} \exp(-\tfrac12 v^{-2} u'A_{n,k}^{-1} u) \dd G(v),
\end{equation}
$u \in \R^k$.  Evaluating \eqref{eq_Fourier_inv_18} at $u = y\bfone_k$ and substituting the result in \eqref{eq_Fourier_inv}, we obtain 
\begin{multline}
\label{eq_E_fn_k}
\E_{\Xi_n} \big[f_{Y_n|\Xi_n}(y)\big]^k \\
= (2\pi)^{-k/2} \, \E \Big[(\det A_{n,k})^{-1/2} \int_0^\infty v^{-k} \exp(- \tfrac12 v^{-2} y^2 \bfone_k' A_{n,k}^{-1} \bfone_k) \dd G(v)\Big].
\end{multline}
Applying in \eqref{eq_E_fn_k} the inequality $\exp(- \tfrac12 v^{-2} y^2 \bfone_k' A_{n,k}^{-1} \bfone_k) \le 1$, $y \in \R$, we find that 
\begin{align*}
\E_{\Xi_n} \big[f_{Y_n|\Xi_n}(y)\big]^k 
&\le (2\pi)^{-k/2} \bigg(\int_0^\infty v^{-k} \dd G(v)\bigg) \E\big[(\det A_{n,k})^{-1/2}\big].
\end{align*}
By H\"older's inequality, 
$$
\int_0^\infty v^{-k} \dd G(v) \le \bigg(\int_0^\infty v^{-d_n} \dd G(v)\bigg)^{k/d_n} < \infty.
$$
Further, since $\E [(\det A_{n_0,k})^{-1/2}] < \infty$ then, by Lemma \ref{lem_An_props}(iii), 
$$
\E [(\det A_{n,k})^{-1/2}] \le \E [(\det A_{n_0,k})^{-1/2}] < \infty
$$
for all $n \ge n_0$.  Therefore the expectation on the right-hand side of \eqref{eq_E_fn_k} converges absolutely and hence, by Tonelli's theorem, the earlier interchanges of expectations and integrals are justified.  

Recall that $\wtilX_{n,1},\ldots,\wtilX_{n,k}$ are continuous and satisfy \ref{condition_1} and \ref{condition_2}.  Therefore as $n \to \infty$, $\wtilX_{n,j}'\wtilX_{n,r} \cip \delta_{j,r} \sigma^2$ for all $j,r=1,\ldots,k$.  Noting that the inverse and the determinant mappings on the cone of positive definite $k \times k$ matrices are continuous functions, it follows that the function 
\begin{align*}
\wtilde{\mcalX}_n &\mapsto (\det A_{n,k})^{-1/2} 
\int_0^\infty v^{-k} \exp(- \tfrac12 v^{-2} y^2 \bfone_k' A_{n,k}^{-1} \bfone_k) \dd G(v) \\
&\equiv (\det \wtilde{\mcalX}_n'\wtilde{\mcalX}_n)^{-1/2} 
\int_0^\infty v^{-k} \exp(- \tfrac12 v^{-2} y^2 \bfone_k' (\wtilde{\mcalX}_n'\wtilde{\mcalX}_n)^{-1} \bfone_k) \dd G(v)
\end{align*}
is continuous since $\wtilde{\mcalX}_n'\wtilde{\mcalX}_n$ is positive definite, almost surely.  By Slutsky's theorem, $A_{n,k} \cip \sigma^2 I_k$; so by the Continuous Mapping Theorem, $\det(A_{n,k}) \cip \sigma^{2k}$ and $\bfone'_k A_{n,k}^{-1} \bfone_k \cip 
k\sigma^{-2}$ as $n \to \infty$.  Applying to \eqref{eq_E_fn_k} the Continuous Mapping Theorem, we find that 
\begin{align}
\label{eq_lim_thm_cond_pdf_k}
\lim_{n \to \infty} \E_{\Xi_n} \big[f_{Y_n|\Xi_n}(y)\big]^k &= (2\pi)^{-k/2} \, \sigma^{-k} \int_0^\infty v^{-k} \exp(- \tfrac12 k \sigma^{-2} y^2 v^{-2}) \dd G(v) \nonumber \\
&\equiv \E_V \big[f_{\mcalN_1(0,\sigma^2 V^2)}(y)\big]^k
\end{align}
for all $y \in \R$, which proves \eqref{eq_conv_cond_pdf_spher} for the case $j=k$.  

Finally, to prove the case in which $j < k$, we apply Lemma \ref{lem_An_props}(iv) to deduce that $\E[(\det A_{n,j})^{-1/2}] < \infty$ for all $j=1,2,\ldots,k$ and all $n \ge n_0$.  Repeating the earlier argument with $k$ replaced by $j$, we obtain \eqref{eq_conv_cond_pdf_spher}.  
$\qed$

\medskip

\noindent
\textit{Proof of Corollary \ref{cor_Lp_conv_cond_pdf}}:  
By Theorem \ref{thm_Lp_conv_cond_pdf_spher}, 
\begin{equation}
\label{eq_lim_cond_pdf_k}
\lim_{n \to \infty} \E_{\Xi_n} \big[f_{Y_n|\Xi_n}(y)\big]^j = \big[f_{\mcalN_1(0,\sigma^2)}(y)\big]^j
\end{equation}
for all $y \in \R$ and all $j=1,\ldots,k$.  Moreover, \eqref{eq_lim_cond_pdf_k} holds trivially for $j=0$.  

Suppose that $k$ is even.  By applying the binomial theorem, we obtain 
\begin{align*}
\E_{\Xi_n} \big|f_{Y_n|\Xi_n}(y) - f_{\mcalN_1(0,\sigma^2)}(y)\big|^k 
&= \sum_{j=0}^k (-1)^j \binom{k}{j} \E_{\Xi_n} \big[f_{Y_n|\Xi_n}(y)\big]^j \big[f_{\mcalN_1(0,\sigma^2)}(y)\big]^{k-j}.
\end{align*}
Letting $n \to \infty$, it follows from \eqref{eq_lim_cond_pdf_k} that 
\begin{equation}
\label{eq_lim_p_even}
\lim_{n \to \infty} \E_{\Xi_n} \big|f_{Y_n|\Xi_n}(y) - f_{\mcalN_1(0,\sigma^2)}(y)\big|^k 
= \big[f_{\mcalN_1(0,\sigma^2)}(y)\big]^k \sum_{j=0}^k (-1)^j \binom{k}{j} = 0.
\end{equation}
By H\"older's inequality, 
\begin{equation}
\label{eq_Ferguson_ineq}
\E_{\Xi_n} \big|f_{Y_n|\Xi_n}(y) - f_{\mcalN_1(0,\sigma^2)}(y)\big|^p \le \big(\E_{\Xi_n} \big|f_{Y_n|\Xi_n}(y) - f_{\mcalN_1(0,\sigma^2)}(y)\big|^{k}\big)^{p/k}.
\end{equation}
Applying \eqref{eq_lim_p_even}, it follows that the left-hand side of \eqref{eq_Ferguson_ineq} converges to $0$ as $n \to \infty$.  This establishes \eqref{eq_Lp_conv_cond_pdf} for the case in which $k$ is even.  

Now suppose that $k$ is odd.  Then by Lemma \ref{lem_An_props}(iv), $\E [(\det A_{n,k-1})^{-1/2}] < \infty$ for all $n \ge n_0$, \textit{i.e.}, the assumptions remain valid with $k$ replaced by $k-1$.  Applying the conclusion obtained for the previous case in which $k$ is even, we deduce that if $k$ is odd then \eqref{eq_Lp_conv_cond_pdf} holds for all $p$ such that $0 < p \le k-1$.  
$\qed$

\medskip

\noindent
\textit{Proof of Theorem \ref{thm_quant_pdf}}: 
It suffices to prove the case in which $j=k$ since all other cases are similar.  
Then, with $V$ and $A_{n,k}$ independent, it follows from \eqref{eq_E_fn_k} that 
$$
\E_{\Xi_n} \big[f_{Y_n|\Xi_n}(y)\big]^k = \E f_{\mcalN_k(0,V^2 A_{n,k})}(y\bfone_k).
$$
By expressing the density $f_{\mcalN_k(0,V^2 A_{n,k})}(\cdot)$ as the inverse Fourier transform of the corresponding characteristic function we obtain, for fixed $V$ and $A_{n,k}$, 
\begin{equation}
\label{eq_inv_FT_device_1}
\E f_{\mcalN_k(0,V^2 A_{n,k})}(y \bfone_k) 
= \E \mcalF_{\goesto{w}{y \bfone_k}}^{-1} \exp(-\tfrac12 V^2 w'A_{n,k}w)
\end{equation}
and, similarly, 
\begin{equation}
\label{eq_inv_FT_device_2}
\E [f_{\mcalN_1(0,\sigma^2 V^2)}(y)]^k \equiv \E f_{\mcalN_k(0,\sigma^2 V^2 I_k)}(y \bfone_k) 
= \E \mcalF_{\goesto{w}{y \bfone_k}}^{-1} \exp(-\tfrac12 \sigma^2 V^2 w'w),
\end{equation}
$y \in \R$.  Subtracting \eqref{eq_inv_FT_device_2} from \eqref{eq_inv_FT_device_1} and  applying the triangle inequality, we obtain  
\begin{align}
\label{eq_quant_diff_pdfs}
\big|\E{}_{\Xi_n} \big[f_{Y_n|\Xi_n}(y)\big]^k & - \E [f_{\mcalN_1(0,\sigma^2 V^2)}(y)]^k\big]\big| \nonumber \\
&\le \E \left|\mcalF_{\goesto{w}{y \bfone_k}}^{-1} \big[\exp(-\tfrac12 V^2 w'A_{n,k}w) - \exp(-\tfrac12 \sigma^2 V^2 w'w)\big]\right|. 
\end{align}

For $\widehat{f} \in L^1(\R^k)$ and $z \in \R^k$, it follows by the integral formula \eqref{eq_FT_inv_integral} for the inverse Fourier transform that 
\begin{equation}
\label{eq_ineq_inv_FT}
\big|\mcalF_{\goesto{w}{z}}^{-1} \widehat{f}(w)\big| 
= \left|(2\pi)^{-k} \int_{\R^k} \exp(-\i z'w) \widehat{f}(w) \dd w\right| 
\le (2\pi)^{-k} \int_{\R^k} |\widehat{f}(w)| \dd w.
\end{equation}
Applying \eqref{eq_ineq_inv_FT} to the right-hand side of \eqref{eq_quant_diff_pdfs}, we obtain an upper bound that does not depend on $y$, and therefore
\begin{multline}
\label{eq_E_diff_pdfs}
\sup_{y \in \R} \big|\E{}_{\Xi_n} \big[f_{Y_n|\Xi_n}(y)\big]^k - \E [f_{\mcalN_1(0,\sigma^2 V^2)}(y)]^k\big| \\
\le (2\pi)^{-k} \, \E \int_{\R^k} \big|\!\exp(-\tfrac12 V^2 w'A_{n,k}w) - \exp(-\tfrac12 \sigma^2 V^2 w'w)\big| \dd w. 
\end{multline}

Denote by $\mathbb{S}^{k \times k}$ the space of real, \textit{symmetric}, $k \times k$ matrices.  For $A = \big(a_{j,l}\big)_{j,l=1}^k \in \mathbb{S}^{k \times k}$, the gradient operator is the $k \times k$ matrix $\nabla = \left(\tfrac12(1+\delta_{j,l}) \partial/\partial a_{j,l}\right)_{j,l=1}^k$, where $\delta_{j,l}$ is Kronecker's delta: $\delta_{j,l} = 1$ if $j=l$ and $\delta_{j,l} = 0$, otherwise.  For $A, A_0 \in \mathbb{S}^{k \times k}$ and a continuously differentiable function $h: \mathbb{S}^{k \times k} \to \R$, the mean-value theorem on $\mathbb{S}^{k \times k}$ \citep[\S 14]{Magnus_Neudecker} states that, for some $\eta \in (0,1)$, 
\begin{equation}
\label{eq_mvt}
h(A) - h(A_0) = \tr \big[(A-A_0) \cdot (\nabla h)\big(\eta A + (1-\eta) A_0\big)\big].
\end{equation}
For fixed $V$ and $w$, set $h(A) = \exp(-\tfrac12 V^2 w'Aw)$, $A \in \mathbb{S}^{k \times k}$; then it is simple to verify that $(\nabla h)(A) = -\tfrac12 V^2 h(A) ww'$.  Setting $A = A_{n,k}$ and $A_0 = \sigma^2 I_k$, applying \eqref{eq_mvt} to the integrand in \eqref{eq_E_diff_pdfs} and simplifying the result, we obtain, for some $\eta \in (0,1)$,  
\begin{multline}
\label{eq_diff_exp}
\exp(-\tfrac12 V^2 w'A_{n,k}w) - \exp(-\tfrac12 \sigma^2 V^2 w'w) \\
= -\tfrac12 V^2 \, w'(A_{n,k} - \sigma^2 I_k)w \cdot \exp\big(\! -\tfrac12 V^2 w'M_1 w\big),
\end{multline}
where $M_1 = \eta A_{n,k} + (1-\eta)\sigma^2 I_k$.  
Inserting \eqref{eq_diff_exp} into \eqref{eq_E_diff_pdfs}, and making the change-of-variables $w \to V^{-1}M_1^{-1/2}w$ in the integral, we obtain   
\begin{align}
\label{eq_diff_exp_2}
\sup_{y \in \R} \big|&\E_{\Xi_n} \big[f_{Y_n|\Xi_n}(y)\big]^k - \E [f_{\mcalN_1(0,\sigma^2 V^2)}(y)]^k\big| \nonumber \\
&\le 2^{-1} (2\pi)^{-k} \, \E(V^{-k}) \, \E\Big[(\det M_1)^{-1/2} \int_{\R^k} |w'M_2w| \cdot \exp\big(-\tfrac12 \|w\|^2\big) \dd w\Big],
\end{align}
where 
\begin{equation}
\label{eq_M2_matrix}
M_2 = M_1^{-1/2} (A_{n,k} - \sigma^2 I_k) M_1^{-1/2}.
\end{equation}
By the Cauchy-Schwarz inequality, $|w'M_2w| \le \|M_2\| \|w\|^2$ for all $w$.  Also, it is a simple Gaussian integral that 
$$
\int_{\R^k} \|w\|^2 \, \exp\big(-\tfrac12 \|w\|^2\big) \dd w =  (2\pi)^{k/2} k.
$$
Therefore it follows from \eqref{eq_diff_exp_2} that 
\begin{multline}
\label{eq_diff_exp_3}
\sup_{y \in \R} \big|\E_{\Xi_n} \big[f_{Y_n|\Xi_n}(y)\big]^k - \E [f_{\mcalN_1(0,\sigma^2 V^2)}(y)]^k\big| \\
\le 2^{-1} (2\pi)^{-k/2} k \, \E(V^{-k}) \, \E\big[(\det M_1)^{-1/2} \|M_2\| \big].
\end{multline}

By \eqref{eq_M2_matrix} and the submultiplicativity property of the Frobenius norm, we have 
\begin{equation}
\label{eq_M_2_norm_bound}
\|M_2\| = \|(A_{n,k} - \sigma^2 I_k) \, M_1^{-1}\|  
\le \|A_{n,k} - \sigma^2 I_k\| \, \|M_1^{-1}\|.
\end{equation}
By \citet[Theorem 8]{Lieb}, for positive definite $A$, the function $A \mapsto \tr(A^{-2})$ is convex.  Therefore 
\begin{align*}
\big\|M_1^{-1}\big\|^2 &= \big\|\big(\eta A_{n,k} + (1-\eta)\sigma^2 I_k\big)^{-1}\big\|^2 \\
&\le \eta \big\|A_{n,k}^{-1}\big\|^2 + (1-\eta) \big\|(\sigma^2 I_k)^{-1}\big\|^2 
\le \max\big\{\big\|A_{n,k}^{-1}\big\|^2,k \sigma^{-4}\big\},
\end{align*}
hence 
$$
\big\|M_1^{-1}\big\| 
\le \max\big\{\big\|A_{n,k}^{-1}\big\|,k^{1/2} \sigma^{-2}\big\}.
$$
Applying the latter inequality to \eqref{eq_M_2_norm_bound}, we obtain 
\begin{equation}
\label{eq_M_2_norm_bound_2}
\big\|M_2\big\| \le \big\|A_{n,k} - \sigma^2 I_k\big\| \cdot \max\big\{\big\|A_{n,k}^{-1}\big\|,k^{1/2} \sigma^{-2}\big\}.
\end{equation}
By \citet[Theorem 7.6.6]{Horn}, for positive definite $A$, the function $A \mapsto \log (\det A)^{-1}$, is convex.  Therefore 
\begin{align*}
\log (\det M_1)^{-1} &\equiv \log \big(\!\det (\eta A_{n,k} + (1-\eta)\sigma^2 I_k)\big)^{-1} \\
&\le \eta \log (\det A_{n,k})^{-1} + (1-\eta) \log (\det \sigma^2 I_k)^{-1} \\
&\le \max\{\log (\det A_{n,k})^{-1},\log \sigma^{-2k}\},
\end{align*}
hence 
\begin{equation}
\label{eq_M1_det_inv}
(\det M_1)^{-1/2} \le \max\{(\det A_{n,k})^{-1/2},\sigma^{-k}\}. 
\end{equation}
Multiplying \eqref{eq_M_2_norm_bound_2} and \eqref{eq_M1_det_inv}, and applying the Cauchy-Schwarz inequality, we obtain 
\begin{align*}
\E&\big[(\det M_1)^{-1/2} \big\|M_2\big\| \big] \\
&\le \big(\E\big\|A_{n,k} - \sigma^2 I_k\big\|^2\big)^{1/2} 
\, \Big(\E\big[\!\max\{(\det A_{n,k})^{-1/2},\sigma^{-k}\} \max\{\big\|A_{n,k}^{-1}\big\|,k^{1/2} \sigma^{-2}\}\big]^2\Big)^{1/2},
\end{align*}
and on substituting this result into \eqref{eq_diff_exp_3} we obtain 
\begin{align}
\label{eq_quant_bound}
\sup_{y \in \R} \big|\E_{\Xi_n} & \big[f_{Y_n|\Xi_n}(y)\big]^k - \E [f_{\mcalN_1(0,\sigma^2 V^2)}(y)]^k\big| \nonumber \\
\le & \ 2^{-1} (2\pi)^{-k/2} k \, \E(V^{-k}) \, \big(\E \big\|A_{n,k} - \sigma^2 I_k\big\|^2\big)^{1/2} \nonumber \\
& \ \cdot 
\Big(\E\big[\max\big\{(\det A_{n,k})^{-1/2},\sigma^{-k}\big\} \cdot \max\big\{\big\|A_{n,k}^{-1}\big\|,k^{1/2} \sigma^{-2}\big\}\big]^{2}\Big)^{1/2}.
\end{align}
Since $A_{n,k} \cip \sigma^2 I_k$ as $n \to \infty$ then 
$(\det A_{n,k})^{-1/2}$ $\cip \sigma^{-k}$ and $\big\|A_{n,k}^{-1}\big\| \cip k^{1/2} \sigma^{-2}$; hence 
$$
\big[\max\{(\det A_{n,k})^{-1/2},\sigma^{-k}\} \cdot \max\big\{\big\|A_{n,k}^{-1}\big\|,k^{1/2} \sigma^{-2}\big\}\big]^2 \cip 
k \sigma^{-2(k+2)}.
$$
Therefore there exists $n_k \in \N$ such that, for all $n \ge n_k$, 
\begin{equation}
\label{eq_quant_bound_2}
\E\big[\max\{(\det A_{n,k})^{-1/2},\sigma^{-k}\} \cdot \max\big\{\big\|A_{n,k}^{-1}\big\|,k^{1/2} \sigma^{-2}\big\}\big]^2 \le 2 k \sigma^{-2(k+2)},
\end{equation}
and then we obtain \eqref{eq_quant_pdf} by applying \eqref{eq_quant_bound_2} to \eqref{eq_quant_bound}.  

To prove \eqref{eq_quant_rate}, we apply the definition of the Frobenius norm to obtain 
\begin{equation}
\label{eq_Ank_diff_sq}
\big\|A_{n,k} - \sigma^2 I_k\big\|^2 = \tr [(A_{n,k} - \sigma^2 I_k)^2] = \tr(A_{n,k}^2 - 2 \sigma^2 A_{n,k} + \sigma^4 I_k).
\end{equation}
By \eqref{eq_A_nk_matrix}, 
$$
\tr(A_{n,k}^2) = \sum_{j=1}^k \sum_{r=1}^k (\wtilX'_{n,j}\wtilX_{n,r})^2 
= \sum_{j=1}^k \|\wtilX_{n,j}\|^4 + \sum_{1 \le j \neq r \le k} (\wtilX'_{n,j}\wtilX_{n,r})^2.
$$
Since $\wtilX_{n,1},\ldots,\wtilX_{n,k}$ are independent copies of $X_n$ then $\E (\|\wtilX_{n,j}\|^4) = \E(\|X_n\|^4)$ for $1 \le j \le r$; also, $\E(\wtilX'_{n,j}\wtilX_{n,r})^2 = \E(X'_n\wtilX_n)^2$ for $1 \le j \neq r \le k$.  Hence 
\begin{equation}
\label{eq_tr_Ank_sq}
\E \tr(A_{n,k}^2) = k \E (\|X_n\|^4) + k(k-1) \E [(X'_n\wtilX_n)^2],
\end{equation}
and, by a similar calculation, $\E \tr(A_{n,k}) = k \E (\|X_n\|^2)$.  Therefore by \eqref{eq_Ank_diff_sq} and \eqref{eq_tr_Ank_sq}, 
$$
\E \big\|A_{n,k} - \sigma^2 I_k\big\|^2 = k \E (\|X_n\|^4) + k(k-1) \E [(X'_n\wtilX_n)^2] - 2 k \sigma^2 \E (\|X_n\|^2) + k \sigma^4,
$$
which reduces to \eqref{eq_quant_rate}.  
$\qed$

\begin{remark}
\label{rem_C_1_prime_sec_4}
{\rm
Suppose that \ref{condition_1} is replaced by \ref{condition_1_prime}.  Let $R_1,\ldots,R_k$ be mutually independent copies of $R$, and assume that $\E(|R|^{-1}) < \infty$.  By Slutsky's and the Continuous Mapping theorems, $A_{n,k} \cip \mcalR_k = \diag(R_1^2,\ldots,R_k^2)$, hence $\bfone'_k A_{n,k}^{-1} \bfone_k \cip 
\sum_{m=1}^k R_m^{-2}$, as $n \to \infty$.  Proceeding similarly to \eqref{eq_E_fn_k}-\eqref{eq_lim_thm_cond_pdf_k} we obtain 
\begin{align}
\label{eq_E_fn_k_cprime}
\E_{\Xi_n} \big[f_{Y_n|\Xi_n}(y)\big]^j 
&\to \E_V \prod_{m=1}^j \E_{R_m} \big[(2\pi)^{-1/2} (V^2 R_m^2)^{-1/2} \exp(-\tfrac12 V^{-2} R_m^{-2} y^2)\big] \nonumber \\
&\equiv \ \E_V \big[\E_R f_{\mcalN_1(0,V^2 R^2)}(y)\big]^j,
\end{align}
for $j=1,\ldots,k$.  Applying the inequality $f_{\mcalN_1(0,V^2 R^2)}(y) \le (2\pi)^{-1/2} (V^2R^2)^{-1/2}$, we have 
\begin{equation}
\label{eq_E_fn_k_cprime_2}
E_V \big[\E_R f_{\mcalN_1(0,V^2 R^2)}(y)\big]^j 
\le (2\pi)^{-j/2} \E(V^{-j}) [\E(|R|^{-1})]^j.
\end{equation}
Since $\E(|R|^{-1}) < \infty$, then \eqref{eq_E_fn_k_cprime_2} is finite, so \eqref{eq_E_fn_k_cprime} extends \eqref{eq_conv_cond_pdf_spher} in Theorem \ref{thm_Lp_conv_cond_pdf_spher}.  

Corollary \ref{cor_Lp_conv_cond_pdf} extends as follows: \textit{Suppose that $\{X_n \in \R^{d_n}, n \ge 1\}$ satisfy \ref{condition_1_prime}, \ref{condition_2}, \ref{condition_3}, and \eqref{eq_integrable_cf_Yn_Xin}, where $\E(|R|^{-1}) < \infty$, and let $\:\Xi_n \sim \mcalN_{d_n}(0,I_{d_n})$.  Then, for all $\:y \in \R$ and all $0 \le p < 2 \lfloor k/2 \rfloor$,}  
\begin{equation}
\label{eq_Lp_conv_cond_pdf_R}
\lim_{n \to \infty} \E_{\Xi_n} \big|f_{Y_n|\Xi_n}(y) - \E f_{\mcalN_1(0,R^2)}(y)\big|^p = 0.
\end{equation}
The proof of \eqref{eq_Lp_conv_cond_pdf_R} is analogous to the proof of \eqref{eq_Lp_conv_cond_pdf}.  

The extension of Theorem \ref{thm_quant_pdf} under \ref{condition_1_prime} is as follows: \textit{Let $X_n$ and $\Xi_n$ satisfy \ref{condition_1_prime}, \ref{condition_2}, \ref{condition_3}, and \eqref{eq_integrable_cf_Yn_Xin}, and assume that $\E(R^{-6}) < \infty$.  Then, for $1 \le j \le k$, there exists $n_j \in \N$ such that, for all $n \ge n_j$, 
$$
\sup_{y \in \R} \big|\E_{\Xi_n} \big[f_{Y_n|\Xi_n}(y)\big]^j - \E_V \big[\E_R f_{\mcalN_1(0,V^2 R^2)}(y)\big]^j\big| \le c_j' \, \big[\E \big\|A_{n,j} - \mcalR_j\big\|^2\big]^{1/2},
$$
where 
\begin{equation}
\label{eq_quant_c_j_R}
c_j' = 2^{-(j+1)/2} \pi^{-j/2} j^{3/2} \E(V^{-j}) (\E R^{-6})^{1/2} (\E R^{-2})^{(k-1)/2}.
\end{equation}
Further, 
$$
\E \big\|A_{n,j} - \mcalR_j\big\|^2 = j \, \E (\|X_n\|^2 - R^2)^2 + j(j-1) [\E (X_n'\wtilX_n)]^2.
$$}
The proof of this extension follows the same arguments in the proof of Theorem \ref{thm_quant_pdf} up to \eqref{eq_quant_bound}, leading to the result that there exists $n_k \in \N$ such that, for all $n \ge n_k$, 
\begin{multline*}
\sup_{y \in \R} \big|\E_{\Xi_n} \big[f_{Y_n|\Xi_n}(y)\big]^k - \E_V [\E_R f_{\mcalN_1(0,V^2 R^2)}(y)]^k\big| \\
\le 2^{-(k+1)/2} \pi^{-k/2} k \, \E(V^{-k}) \, \big(\E \big\|A_{n,k} - \mcalR_k\big\|^2\big)^{1/2} \cdot 
\big(\E \big[(\det \mcalR_k)^{-1} \big\|\mcalR_k^{-1}\big\|^2\big]\big)^{1/2}.
\end{multline*}
Also, the constant $c_j'$ in \eqref{eq_quant_c_j_R} is finite since 
$$
\E \big[(\det \mcalR_k)^{-1} \big\|\mcalR_k^{-1}\big\|^2\big] = \E \big[(R_1 \cdots R_k)^{-2} \, (R_1^{-4} + \cdots + R_k^{-4})\big] = k (\E R^{-6}) (\E R^{-2})^{k-1}.
$$
}\end{remark}

\section{Properties of the cumulative distribution function of \texorpdfstring{$\boldsymbol{Y_n|\Xi_n}$}{Yn}}
\label{sec_Lp_conv_cond_cdf}

In this section we obtain conditions under which $F_{Y_n|\Xi_n}$, the distribution function of $Y_n|\Xi_n$, converges uniformly to the distribution function of a mixture of normal distributions.  This result is motivated by classical statistical inference, in which the well-known Glivenko-Cantelli theorem establishes the uniform convergence of an empirical distribution function to its population counterpart.  

In the ensuing results, we retain the notation of Section \ref{sec_Lp_conv_cond_pdf}; thus $G$ denotes the distribution function defined by \eqref{eq_Schoenberg} and $V$ denotes the random variable corresponding to $G$.  Further we denote by $F_{\mcalN_1(0,\sigma^2)}$ the distribution function of the $\mcalN_1(0,\sigma^2)$ distribution.

\subsection{Convergence properties of the cumulative distribution function}
\label{subsec_convergence_props_cdf}

\begin{theorem}
\label{thm_Lp_conv_cond_cdf_spher}
Let $\{X_n \in \R^{d_n}, n \ge 1\}$ be continuous random vectors that satisfy \ref{condition_1}, \ref{condition_2}, \ref{condition_3}, and \eqref{eq_integrable_cf_Yn_Xin}.  Let $\{\Xi_n \in \R^{d_n}, n \ge 1\}$ be spherically symmetric modulating random vectors that satisfy \eqref{eq_Xi_n_cf} and \eqref{eq_G_dn_integrable} and are independent of $\{X_n, n \ge 1\}$, and let $Y_n = \Xi_n'X_n$, $n \ge 1$.  Then, for all $j=1,\ldots,k$, 
\begin{equation}
\label{eq_Lp_conv_cond_cdf_spher}
\lim_{n \to \infty} \sup_{y \in \R} \Big|\E_{\Xi_n} \big[F_{Y_n|\Xi_n}(y)\big]^j - \E_V \big[F_{\mcalN_1(0,\sigma^2 V^2)}(y)\big]^j \Big| = 0.
\end{equation}
\end{theorem}

For the case $p=2$, \citet[Theorem 2.6]{Bagyan} established the following result.  

\begin{corollary}
\label{cor_Lp_conv_cond_cdf_gaussian}
Suppose that the continuous random vectors $\{X_n \in \R^{d_n}, n \ge 1\}$ satisfy \ref{condition_1}, \ref{condition_2}, \ref{condition_3}, and \eqref{eq_integrable_cf_Yn_Xin}.  Let $\Xi_n \sim \mcalN_{d_n}(0,I_{d_n})$, $n \ge 1$, and suppose that $\{\Xi_n, n \ge 1\}$ and $\{X_n, n \ge 1\}$ are independent.  Then for all $0 \le p < 2 \lfloor k/2 \rfloor$, 
\begin{equation}
\label{eq_Lp_conv_cond_cdf}
\lim_{n \to \infty} \, \sup_{y \in \R} \E_{\Xi_n} \big|F_{Y_n|\Xi_n}(y) - F_{\mcalN_1(0,\sigma^2)}(y)\big|^p = 0.
\end{equation}
\end{corollary}

Recalling the well-known result that the L\'evy metric is dominated by the supremum (\textit{i.e.}, Kolmogorov) metric, it follows that Corollary \ref{cor_Lp_conv_cond_cdf_gaussian} remains valid if distances between distribution functions are measured using the L\'evy metric.  

In the next result, we provide a quantitative version of Theorem \ref{thm_Lp_conv_cond_cdf_spher}.  It is also evident that this result represents a Lipschitz continuity property of $F_{Y_n|\Xi_n}(\cdot)$.  

\begin{theorem}
\label{thm_quant_cdf}
Suppose that $X_n$ and $\Xi_n$ satisfy the assumptions of Theorem \ref{thm_Lp_conv_cond_cdf_spher} and that $1 \le j \le k$.  Let $y, a \in \R$ and let $c_j$ be the constant defined in \eqref{eq_quant_c_k}.  Then there exists $n_j \in \N$ such that, for all $n \ge n_j$, 
\begin{multline}
\label{eq_quant_cdf_lipschitz}
\big|\E_{\Xi_n} \big[F_{Y_n|\Xi_n}(y) - F_{Y_n|\Xi_n}(a)\big]^j - \E_V \big[F_{\mcalN_1(0,\sigma^2 V^2)}(y) - F_{\mcalN_1(0,\sigma^2 V^2)}(a)\big]^j\big| \\
\le c_j \, |y-a|^j \, \big[\E \big\|A_{n,j} - \sigma^2 I_j\big\|^2\big]^{1/2}.
\end{multline}
\end{theorem}

\subsection{Proofs}
\label{subsec_cdf_proofs}

\noindent
\textit{Proof of Theorem \ref{thm_Lp_conv_cond_cdf_spher}}:  
Since $\Xi_n$ is independent of $X_n$ then the conditional characteristic function of $Y_n$ given $\Xi_n$ is 
$$
\varphi_{Y_n|\Xi_n}(t) = \E_{Y_n|\Xi_n} \exp(\i tY_n) = \E_{X_n|\Xi_n} \exp(\i t\Xi_n'X_n),
$$
$t \in \R$.   
Therefore 
\begin{align*}
F_{Y_n|\Xi_n}(y) 
&= \int_{-\infty}^y f_{Y_n|\Xi_n}(z) \dd z \\
&= \int_{-\infty}^y (\mcalF^{-1}\varphi_{Y_n|\Xi_n})(z) \dd z 
= \int_{-\infty}^y \big[\mcalF_{\goesto{t}{z}}^{-1} \E_{X_n|\Xi_n}\exp(\i  t\Xi_n'X_n)\big] \dd z.
\end{align*}
Let $\wtilX_{n,1},\ldots,\wtilX_{n,k}$ be independent copies of $X_n$; then 
$$
\big[F_{Y_n|\Xi_n}(y)\big]^k = \int_{-\infty}^y \cdots \int_{-\infty}^y \prod_{j=1}^k \mcalF_{\goesto{t_j}{z_j}}^{-1} \E_{\wtilX_{n,j}|\Xi_n} \exp(\i  t_j\Xi_n'\wtilX_{n,j}) \dd z_j.
$$
Set $u = (z_1,\ldots,z_k)'$, and formally interchange $\E_{\Xi_n}$ with the multiple integral and the operators $\mcalF_{\goesto{t_j}{z_j}}^{-1} \E_{\wtilX_{n,j}}$, $j=1,\ldots,k$; then we obtain 
$$
\E_{\Xi_n} \big[F_{Y_n|\Xi_n}(y)\big]^k = \idotsint\limits_{u \in (-\infty,y]^k} \Big(\prod_{j=1}^k \mcalF_{\goesto{t_j}{z_j}}^{-1} \E_{\wtilX_{n,j}}\Big) \E_{\Xi_n} \exp\Big(\i \sum_{j=1}^k t_j\Xi_n'\wtilX_{n,j}\Big) \dd u.
$$
Since $\Xi_n$ is spherically symmetric then it follows from \eqref{eq_Ank_psi} that 
$$
\E_{\Xi_n} \big[F_{Y_n|\Xi_n}(y)\big]^k = \idotsint\limits_{u \in (-\infty,y]^k} \Big(\prod_{j=1}^k \mcalF_{\goesto{t_j}{z_j}}^{-1} \E_{\wtilX_{n,j}}\Big) \psi(w'A_{n,k}w)  \dd u,
$$
and by formally interchanging Fourier transforms and expectations we obtain 
\begin{equation}
\label{eq_E_power_cdf_spherical_1}
\E_{\Xi_n} \big[F_{Y_n|\Xi_n}(y)\big]^k = \idotsint\limits_{u \in (-\infty,y]^k} \E_{A_{n,k}} \mcalF_{\goesto{t_1}{z_1}}^{-1} \cdots \mcalF_{\goesto{t_k}{z_k}}^{-1} \psi(w'A_{n,k}w)  \dd u.
\end{equation}
For $d_n \ge k$, it follows as in \eqref{eq_Fourier_inv_18} that 
\begin{multline*}
\mcalF_{\goesto{t_1}{z_1}}^{-1} \cdots \mcalF_{\goesto{t_k}{z_k}}^{-1} \psi(w' A_{n,k} w) \\
= (2\pi)^{-k/2} (\det A_{n,k})^{-1/2} \int_0^\infty v^{-k} \exp(-\tfrac12 v^{-2} u'A_{n,k}^{-1} u) \dd G(v).
\end{multline*}
Substituting this result in \eqref{eq_E_power_cdf_spherical_1} and interchanging expectations and integrals, we obtain 
\begin{align}
\label{eq_E_power_cdf_spherical_2}
\E&{}_{\Xi_n} \big[F_{Y_n|\Xi_n}(y)\big]^k \nonumber \\
&= \E_{A_{n,k}} \int_0^\infty \idotsint\limits_{u \in (-\infty,y]^k} (2\pi)^{-k/2} (\det A_{n,k})^{-1/2} v^{-k} \exp(-\tfrac12 v^{-2} u'A_{n,k}^{-1} u) \dd u \dd G(v).
\end{align}
Let $(Z_1,\ldots,Z_k)'|A_{n,k} \sim \mcalN_k(0,v^2 A_{n,k})$; then by \eqref{eq_E_power_cdf_spherical_2}, 
\begin{align*}
\E_{\Xi_n} \big[F_{Y_n|\Xi_n}(y)\big]^k &= 
\E_{A_{n,k}} \int_0^\infty \P(Z_1 \le z,\ldots,Z_k \le z | A_{n,k}) \dd G(v) \\
&\le \E \int_0^\infty \dd G(v) = \E(1) = 1,
\end{align*}
proving that \eqref{eq_E_power_cdf_spherical_2} converges.  Therefore by the Fubini-Tonelli theorem, all the foregoing interchanges of integrals and expectations are justified by the absolute convergence of the resulting integral.  

Since the integrand in \eqref{eq_E_power_cdf_spherical_2} is bounded and continuous then, by the Dominated Convergence theorem and the property $A_{n,k} \cip \sigma^2 I_k$, it follows that, as $n \to \infty$, 
\begin{eqnarray*}
\E_{\Xi_n} \big[F_{Y_n|\Xi_n}(y)\big]^k 
& \to & \int_0^\infty \idotsint\limits_{u \in (-\infty,y]^k} (2\pi)^{-k/2} \sigma^{-k} v^{-k} \exp(-\tfrac12 v^{-2} \sigma^{-2} u'u) \dd u \dd G(v) \\
& = & \E_V \big[F_{\mcalN_1(0,\sigma^2 V^2)}(y)\big]^k.
\end{eqnarray*}

Next, we follow the last part of the proof of Theorem \ref{thm_Lp_conv_cond_pdf_spher}.  Starting with the assumption that $\E [(\det A_{n_0,k})^{1/2}] < \infty$ for some $n_0$, we apply Lemma \ref{lem_An_props}(iii,iv) to deduce that $\E [(\det A_{n,j})^{1/2}] < \infty$ for all $n \ge n_0$ and all $j=1,\ldots,k-1$.  By repeating the earlier arguments with $k$ replaced successively by $k-1$, we obtain 
\begin{equation}
\label{eq_lim_k_power_cdf_spherical} 
\E_{\Xi_n} \big[F_{Y_n|\Xi_n}(y)\big]^j = \E \big[F_{\mcalN_1(0,\sigma^2 V^2)}(y)\big]^j
\end{equation}
for all $y \in \R$ and all $j=k,k-1,k-2\ldots,1$.  

To show that the convergence in \eqref{eq_lim_k_power_cdf_spherical} is uniform in $y$ we note that the function $\E_{\Xi_n} \big[F_{Y_n|\Xi_n}(y)\big]^j$, $y \in \R$, also is a cumulative distribution function.  Indeed, since $\Xi_n$ is independent of $\wtilX_{n,1},\ldots,\wtilX_{n,k}$ then 
\begin{align*}
\E_{\Xi_n} \big[F_{Y_n|\Xi_n}(y)\big]^j &= \E_{\Xi_n} \big[\P(Y_n \le y|\Xi_n)\big]^j \nonumber \\
&= \E_{\Xi_n} \P\big(\Xi_n'\wtilX_{n,1} \le y,\ldots,\Xi_n'\wtilX_{n,j} \le y|\Xi_n\big) \nonumber \\
&= \P\big(\max\{\Xi_n'\wtilX_{n,1},\ldots,\Xi_n'\wtilX_{n,j}\} \le y\big),
\end{align*}
which clearly is a cumulative distribution function.  Consequently, $\E_{\Xi_n} \big[F_{Y_n|\Xi_n}(y)\big]^j$ converges to $0$ as $y \to -\infty$ and to $1$ as $y \to \infty$; and by a similar argument, it is also evident that $\E_V \big[F_{\mcalN_1(0,\sigma^2 V^2)}(y)\big]^j$, $y \in \R$, is a cumulative distribution function, and it converges to $0$ as $y \to -\infty$ and to $1$ as $y \to \infty$.  

As the distribution function $\E_{\Xi_n} \big[F_{Y_n|\Xi_n}(y)\big]^j$ converges pointwise to the distribution function $\E_V \big[F_{\mcalN_1(0,\sigma^2 V^2)}(y)\big]^j$, and since both functions attain the same values as $y \to \pm\infty$ then, by \citet[p.~338, Theorem 9.1.6]{Kawata}, $\E_{\Xi_n} \big[F_{Y_n|\Xi_n}(y)\big]^j$ converges uniformly to $\E_V \big[F_{\mcalN_1(0,\sigma^2 V^2)}(y)\big]^j$ as $n \to \infty$.  Therefore \eqref{eq_Lp_conv_cond_cdf_spher} is established.  
$\qed$

\medskip

\noindent
\textit{Proof of Corollary \ref{cor_Lp_conv_cond_cdf_gaussian}}:  
Since $\Xi_n \sim \mcalN_{d_n}(0,I_{d_n})$ then $G$ is singular, with $V = 1$, almost surely, and \eqref{eq_G_dn_integrable} holds trivially.  Therefore, by Theorem \ref{thm_Lp_conv_cond_cdf_spher}, 
\begin{equation}
\label{eq_unif_lim_cond_pdf_k}
\sup_{y \in \R} \left|\E_{\Xi_n} \big[F_{Y_n|\Xi_n}(y)\big]^j - \big[F_{\mcalN_1(0,\sigma^2)}(y)\big]^j\right| \to 0
\end{equation}
as $n \to \infty$, for all $j=1,\ldots,k$.  Also, \eqref{eq_unif_lim_cond_pdf_k} obviously holds for $j=0$.  

Suppose that $k$ is even.  Applying the binomial theorem, we obtain 
\begin{align*}
\E_{\Xi_n} \big|F_{Y_n|\Xi_n}(y) - F_{\mcalN_1(0,\sigma^2)}(y)\big|^k &\equiv \E_{\Xi_n} \big[F_{Y_n|\Xi_n}(y) - F_{\mcalN_1(0,\sigma^2)}(y)\big]^k \\
&= \sum_{j=0}^k (-1)^j \binom{k}{j} \E_{\Xi_n} \big[F_{Y_n|\Xi_n}(y)\big]^j \big[F_{\mcalN_1(0,\sigma^2)}(y)\big]^{k-j}.
\end{align*}
Since $E_{\Xi_n} \big[F_{Y_n|\Xi_n}(y)\big]^j$ converges uniformly to $F_{\mcalN_1(0,\sigma^2)}(y)$ as $n \to \infty$ then, by \eqref{eq_unif_lim_cond_pdf_k}, 
\begin{equation}
\label{eq_lim_p_even_2}
\E_{\Xi_n} \big|F_{Y_n|\Xi_n}(y) - F_{\mcalN_1(0,\sigma^2)}(y)\big|^k \to \big[F_{\mcalN_1(0,\sigma^2)}(y)\big]^k \sum_{j=0}^k (-1)^j \binom{k}{j} \equiv 0,
\end{equation}
with uniform convergence in $y$.  By H\"older's inequality, for $p \le k$, 
\begin{equation}
\label{eq_Ferguson_ineq_2}
\E_{\Xi_n} \big|F_{Y_n|\Xi_n}(y) - F_{\mcalN_1(0,\sigma^2)}(y)\big|^p \le \big(\E_{\Xi_n} \big|F_{Y_n|\Xi_n}(y) - F_{\mcalN_1(0,\sigma^2)}(y)\big|^{k}\big)^{p/k}.
\end{equation}
Applying \eqref{eq_lim_p_even_2}, it follows that the left-hand side of \eqref{eq_Ferguson_ineq_2} converges uniformly to $0$ as $n \to \infty$.  This establishes \eqref{eq_Lp_conv_cond_cdf} for the case in which $k$ is even.  

Finally, for $k$ odd, we proceed as before, applying Lemma \ref{lem_An_props}(iii,iv) to reduce the argument to the case in which $k$ is replaced by $k-1$.
$\qed$

\medskip

\noindent
\textit{Proof of Theorem \ref{thm_quant_cdf}}:  Assume, without loss of generality, that $y \ge a$.  Then, proceeding analogously to the proof of \eqref{eq_E_power_cdf_spherical_2}, one deduces that 
\begin{align}
\label{eq_cdf_lipschitz_1}
\E_{\Xi_n} \big[F_{Y_n|\Xi_n}(y) - F_{Y_n|\Xi_n}(a)\big]^k &= \E \idotsint\limits_{u \in (a,y]^k} f_{\mcalN_k(0,V^2 A_{n,k})}(u) \dd u \nonumber \\
&= \E \idotsint\limits_{u \in (a,y]^k} \mcalF_{\goesto{w}{u}}^{-1} \exp(-\tfrac12 V^2 w'A_{n,k}w) \dd u,
\end{align}
where the second equality is obtained by applying the inverse Fourier transform device from \eqref{eq_inv_FT_device_1}.  Proceeding similarly when $A_{n,k}$ is replaced by $\sigma^2 I_k$, and using \eqref{eq_inv_FT_device_2}, we obtain 
\begin{align}
\label{eq_cdf_lipschitz_2}
\E \big[F_{\mcalN_1(0,\sigma^2 V^2)}(y) - F_{\mcalN_1(0,\sigma^2 V^2)}(a)\big]^k &= \E \idotsint\limits_{u \in (a,y]^k} f_{\mcalN_k(0,\sigma^2 V^2 I_k)}(u) \dd u \nonumber \\
&= \E \idotsint\limits_{u \in (a,y]^k} \mcalF_{\goesto{w}{u}}^{-1} \exp(-\tfrac12 \sigma^2 V^2 w'w) \dd u.
\end{align}
Subtracting \eqref{eq_cdf_lipschitz_2} from \eqref{eq_cdf_lipschitz_1} and then applying the triangle inequality, we obtain 
\begin{align}
\label{eq_cdf_lipschitz_3}
\Big|\E_{\Xi_n} &\big[F_{Y_n|\Xi_n}(y) - F_{Y_n|\Xi_n}(a)\big]^k - \E [F_{\mcalN_1(0,\sigma^2 V^2)}(y) - F_{\mcalN_1(0,\sigma^2 V^2)}(a)\big]^k\Big| \nonumber \\
&\le \E \idotsint\limits_{u \in (a,y]^k} \big|\mcalF_{\goesto{w}{u}}^{-1} \big[\exp(-\tfrac12 V^2 w'A_{n,k}w) - \exp(-\tfrac12 \sigma^2 V^2 w'w)\big]\big| \dd u.
\end{align}
Now applying \eqref{eq_ineq_inv_FT}, it follows that \eqref{eq_cdf_lipschitz_3} is bounded above by  
\begin{align*}
(2\pi)^{-k} \, & \E \idotsint\limits_{u \in (a,y]^k} \int_{\R^k} \big|\exp(-\tfrac12 V^2 w'A_{n,k}w) - \exp(-\tfrac12 \sigma^2 V^2 w'w)\big| \dd w \dd u \\
&= (2\pi)^{-k} (y-a)^k \, \E \int_{\R^k} \big|\exp(-\tfrac12 V^2 w'A_{n,k}w) - \exp(-\tfrac12 \sigma^2 V^2 w'w)\big| \dd w.
\end{align*}
Noting that the latter expectation is precisely the expectation on the right-hand side of \eqref{eq_E_diff_pdfs}, we apply the upper bound obtained in \eqref{eq_quant_bound} for that expectation, \textit{viz.}, 
\begin{align*}
\E \int_{\R^k} & \big|\exp(-\tfrac12 V^2 w'A_{n,k}w) - \exp(-\tfrac12 \sigma^2 V^2 w'w)\big| \dd w \\
&\le (2\pi)^{k/2} k^{5/4} \E(V^{-k}) \, \big(\E\|A_{n,k} - \sigma^2 I_k\|^2\big)^{1/2} \\
& \quad \cdot 
\Big(\E\big[\max\{(\det A_{n,k})^{-1/2},\sigma^{-k}\} \cdot \max\{\|A_{n,k}^{-1/2}\|,k^{1/2} \sigma^{-1}\}\big]^2\Big)^{1/2},
\end{align*}
and also apply \eqref{eq_quant_bound_2}; then we obtain \eqref{eq_quant_cdf_lipschitz}.  
$\qed$

\begin{remark}
\label{rem_C_1_prime_sec_5}
{\rm
Suppose that \ref{condition_1} is replaced by \ref{condition_1_prime}.  Then the extension of Theorem \ref{thm_Lp_conv_cond_cdf_spher} is as follows: \textit{Let $\{X_n \in \R^{d_n}, n \ge 1\}$ be continuous random vectors that satisfy \ref{condition_1_prime}, \ref{condition_2}, \ref{condition_3}, and \eqref{eq_integrable_cf_Yn_Xin}.  Let $\{\Xi_n \in \R^{d_n}, n \ge 1\}$ be spherically symmetric modulating random vectors that satisfy \eqref{eq_Xi_n_cf} and \eqref{eq_G_dn_integrable} and are independent of $\{X_n, n \ge 1\}$, and define $Y_n = \Xi_n'X_n$, $n \ge 1$.  Then, for all $j=1,\ldots,k$, 
$$
\lim_{n \to \infty} \sup_{y \in \R} \Big|\E_{\Xi_n} \big[F_{Y_n|\Xi_n}(y)\big]^j - \E_V \big[\E_R F_{\mcalN_1(0,V^2 R^2)}(y)\big]^j \Big| = 0.
$$}
The proof of the extension is entirely similar to the proof of Theorem \ref{thm_Lp_conv_cond_cdf_spher}, where we use the same methods as in Remark \ref{rem_C_1_prime_sec_4}.  

Corollary \ref{cor_Lp_conv_cond_cdf_gaussian} extends under \ref{condition_1_prime} as follows: \textit{Suppose that the continuous random vectors $\{X_n \in \R^{d_n}, n \ge 1\}$ satisfy \ref{condition_1_prime}, \ref{condition_2}, \ref{condition_3}, and \eqref{eq_integrable_cf_Yn_Xin}.  Suppose that $\Xi_n \sim \mcalN_{d_n}(0,I_{d_n})$, $n \ge 1$, and that $\{\Xi_n, n \ge 1\}$ and $\{X_n, n \ge 1\}$ are independent.  Then for all $0 \le p < 2 \lfloor k/2 \rfloor$,} 
$$
\lim_{n \to \infty} \, \sup_{y \in \R} \E_{\Xi_n} \big|F_{Y_n|\Xi_n}(y) - \E_R F_{\mcalN_1(0,R^2)}(y)\big|^p = 0.
$$

Theorem \ref{thm_quant_cdf} extends under \ref{condition_1_prime} as follows: 
\textit{Let $\{X_n \in \R^{d_n}, n \ge 1\}$ and $\{\Xi_n \in \R^{d_n}, n \ge 1\}$ be mutually independent random vectors that satisfy \ref{condition_1_prime}, \ref{condition_2}, \ref{condition_3}, and \eqref{eq_integrable_cf_Yn_Xin}, and assume that $\E(R^{-6}) < \infty$.  Let $y, a \in \R$ and let $c_j'$ be the constant defined in \eqref{eq_quant_c_j_R}.  Then for $1 \le j \le k$, there exists $n_j \in \N$ such that, for all $n \ge n_j$,} 
\begin{multline*}
\big|\E_{\Xi_n} \big[F_{Y_n|\Xi_n}(y) - F_{Y_n|\Xi_n}(a)\big]^j - \E_V \big[\E_R \big(F_{\mcalN_1(0,V^2 R^2)}(y) - F_{\mcalN_1(0,V^2 R^2)}(a)\big)\big]^j\big| \\
\le c_j' \, |y-a|^j \, \big[\E \big\|A_{n,j} - \mcalR_j\big\|^2\big]^{1/2}.
\end{multline*}
}\end{remark}

\section{Examples of distributions satisfying \ref{condition_3} and \texorpdfstring{{\color{blue}{(}}\ref{eq_integrable_cf_Yn_Xin}{\color{blue}{)}}}{41}}
\label{sec_examples_2}

We now show that the distributions considered in Examples \ref{ex_Xn_Bingham}-\ref{ex_Xn_gaussian} satisfy the assumptions in Theorems \ref{thm_Lp_conv_cond_pdf_spher} and \ref{thm_Lp_conv_cond_cdf_spher}.  Since we have already verified \ref{condition_1} and \ref{condition_2} for those examples then we need only to verify the integrability requirements \ref{condition_3} and \eqref{eq_integrable_cf_Yn_Xin}.  Further, we provide examples of vectors $\Xi_n$ that satisfy the preceding results.

\begin{example} {\rm 
\label{ex_Xn_gaussian_3}
(Continuation of Example \ref{ex_Xn_gaussian}): 
Let $X_n \sim \mcalN_{d_n}(0,\Sigma_n)$ where $\Sigma_n$ is positive definite.  Then $\varphi_{X_n}(t\theta) = \exp(-t^2\theta'\Sigma_n\theta/2)$, $t \in \R$ and $\theta \in \mcalS^{d_n-1}$.  So \eqref{eq_integrable_cf_Yn_Xin} holds trivially, and $\varphi_{Y_n|\Xi_n}$ is integrable for almost all values of $\Xi_n$.  

Let $H_n$ be a $d_n \times d_n$ orthogonal matrix such that $H_n\Sigma_n H_n'$ is diagonal.  Since $A_{n,k}$ is unchanged when each $\wtilX_{n,j}$ is transformed to $H_n\wtilX_{n,j}$, $j=1,\ldots,k$ then, without loss of generality, we assume that $\Sigma_n$ is diagonal and denote by $\lambda_{n;1},\ldots,\lambda_{n;d_n}$ its diagonal entries.  
Since $\Sigma_n$ is diagonal then all $kd_n$ entries, $\{\wtilX_{n,j;m}, 1 \le j \le k, 1 \le m \le d_n\}$, of the matrix $\wtilde{\mcalX}_n$ are mutually independent.  Also $(\wtilX_{n,1;m},\ldots,\wtilX_{n,k;m})' \sim \mcalN_k(0,\lambda_{n;m} I_k)$, $1 \le m \le d_n$, so the $k \times k$ matrix 
$$
W_{n,m} := \lambda_{n;m}^{-1} (\wtilX_{n,1;m},\ldots,\wtilX_{n,k;m})'(\wtilX_{n,1;m},\ldots,\wtilX_{n,k;m})
$$
has a Wishart distribution with $1$ degree-of-freedom and matrix parameter $I_k$, written $W_{n,m} \sim \mcalW_k(1,I_k)$.  Moreover, $W_{n,1},\ldots,W_{n,d_n}$ are mutually independent and identically distributed as $\mcalW_k(1,I_k)$ and, by \eqref{eq_An_decomp}, 
\begin{equation}
\label{eq_An_decomp_2}
A_{n,k} = \sum_{m=1}^{d_n} \lambda_{n;m} W_{n,m}.
\end{equation}

Let $\lambda_{n;0} = \min\{\lambda_{n;m}, 1 \le m \le d_n\}$, the smallest eigenvalue of $\Sigma_n$, and define $W_n = \sum_{m=1}^{d_n} W_{n,m}$.  Then $W_n \sim \mcalW_k(d_n,I_k)$, which is a nonsingular Wishart distribution since $d_n \ge k$, and by \eqref{eq_An_decomp_2}, 
$$
A_{n,k} = \sum_{m=1}^{d_n} \lambda_{n;0} W_{n,m} + \sum_{m=1}^{d_n} (\lambda_{n;m} - \lambda_{n;0}) W_{n,m} = \lambda_{n;0} W_n + \sum_{m=1}^{d_n} (\lambda_{n;m} - \lambda_{n;0})W_{n,m},
$$
a nonnegative linear combination of positive semidefinite matrices.  Therefore 
$$
\det(A_{n,k}) 
\ge \det\big(\lambda_{n;0} W_n\big) = \lambda_{n;0}^k \det(W_n).
$$
By \citet[p.~101]{Muirhead} we have, for $d_n \ge k+1$, 
\begin{equation}
\label{eq_det_B}
\E[(\det A_{n,k})^{-1/2}] \le \lambda_{n;0}^{-k/2} \, \E [(\det W_n)^{-1/2}] 
= 2^{-k/2} \lambda_{n;0}^{-k/2} \prod_{j=1}^k \frac{\Gamma(\frac12(d_n-j))}{\Gamma(\frac12(d_n-j+1))},
\end{equation}
so $\E[(\det A_{n,k})^{-1/2}] < \infty$ for all $n$ such that $d_n \ge k+1$, and hence \ref{condition_3} holds.  

Note that for the case in which $\Sigma_n = d_n^{-1} \sigma^2 I_{d_n}$, which is the special case of \eqref{eq_Xn_gaussian_3} with $r=0$, it follows from \eqref{eq_An_decomp_2} that $A_{n,k} \sim \mcalW_k(d_n,d_n^{-1} \sigma^2 I_k)$.  Then by \eqref{eq_det_B}, 
\begin{equation}
\label{eq_det_An_general_k}
\E[(\det A_{n,k})^{-1/2}] 
= 2^{-k/2} d_n^{k/2} \sigma^{-k} \prod_{j=1}^k \frac{\Gamma(\frac12(d_n-j))}{\Gamma(\frac12(d_n-j+1))}.
\end{equation}
By applying Stirling's approximation for the gamma function, it follows from \eqref{eq_det_An_general_k} that $\E[(\det A_{n,k})^{-1/2}] \to \sigma^{-k}$ as $n \to \infty$.  This result is consistent with \ref{condition_1} and \ref{condition_2} since, under those assumptions, $A_{n,k} \cip \sigma^2 I_k$ and therefore $\E[(\det A_{n,k})^{-1/2}] \to \det(\sigma^2 I_k)^{-1/2} \equiv \sigma^{-k}$ as $n \to \infty$.  
}\end{example}

\begin{example}{\rm
\label{ex_Bingham_contd_2}
(Continuation of Example \ref{ex_Xn_Bingham}):  Let $\Theta_n$ be Bingham-distributed with matrix parameter $\Sigma_n$.  As noted earlier, the density function \eqref{eq_Bingham_pdf} remains unchanged if $\Sigma_n$ is replaced by $\Sigma_n - \tau I_{d_n}$, for any constant $\tau \in \R$.  By choosing $\tau$ suitably large we may assume, without loss of generality, that $\Sigma_n$ is negative definite; and now we define $\Lambda_n = (-2\Sigma_n)^{-1}$, so that $\Lambda_n$ is positive definite.  

As noted by \citet[p.~841]{Bingham85} and \citet{Kume_Walker}, the Bingham distribution arises by constraining the multivariate normal distribution to $\mcalS^{d_n-1}$; \textit{i.e.}, if $Z_n \sim \mcalN_{d_n}(0,\Lambda_n)$ then $\Theta_n \eqdist Z_n{\big|}\{\|Z_n\|=1\}$.  Therefore for $t \in \R$ and $\theta \in \mcalS^{d_n-1}$, 
$$
\varphi_{\Theta_n}(t\theta) = \E \exp(\i t\theta'\Theta_n) 
= \E_{Z_n|\{\|Z_n\|=1\}} \exp(\i t\theta'Z_n).
$$
For fixed $\theta \in \mcalS^{d_n-1}$, suppose that $\int_{-\infty}^\infty |\varphi_{\Theta_n}(t\theta)| \dd t$ diverges.  By the change-of-variable $t \to st$, where $s > 0$, it follows that $\int_{-\infty}^\infty |\varphi_{\Theta_n}(st\theta)| \dd t$ diverges for all $s$.  Note that 
\begin{align*}
\int_{-\infty}^\infty |\varphi_{\Theta_n}(st\theta)| \dd t &= \int_{-\infty}^\infty |\E_{Z_n|\{\|Z_n\|=1\}} \exp(\i st\theta'Z_n)| \dd t \\
&= \int_{-\infty}^\infty |\E_{Z_n|\{\|Z_n\|=s\}} \exp(\i t\theta'Z_n)| \dd t,
\end{align*}
and then integrating with respect to $s$, we deduce that 
$$
\int_{-\infty}^\infty |\E \exp(\i t\theta'Z_n | \{\|Z_n\| \le s\})| \dd t
$$
diverges for all $s > 0$.  Now letting $s \to \infty$, it follows that 
\begin{equation}
\label{eq_Yn_cf_integral}
\int_{-\infty}^\infty |\E \exp(\i t\theta'Z_n)| \dd t
\end{equation}
also diverges.  However since $Z_n \sim \mcalN_{d_n}(0,\Lambda_n)$ then \eqref{eq_Yn_cf_integral} converges for all $\theta \in \mcalS^{d_n-1}$, as shown in Example \ref{ex_Xn_gaussian_3}.  Therefore we deduce, by contradiction, that \eqref{eq_integrable_cf_Yn_Xin} holds for the Bingham distributions.  

Let $\wtilde{\Theta}_{n,1},\ldots,\wtilde{\Theta}_{n,k}$ and $\wtilde{Z}_{n,1},\ldots,\wtilde{Z}_{n,k}$ be independent copies of $\Theta_n$ and $Z_n$, respectively, and define the $k \times k$ matrices $B_{n,k} = (\wtilde{\Theta}_{n,j}'\wtilde{\Theta}_{n,r})_{j,r=1}^k$ and $C_{n,k} = (\wtilde{Z}_{n,j}'\wtilde{Z}_{n,r})_{j,r=1}^k$.  Again using the relationship between the Bingham and the multivariate normal distributions, we obtain 
$$
\E[(\det B_{n,k})^{-1/2}] = \E\big[(\det C_{n,k})^{-1/2} \big| \{\|\wtilde{Z}_{n,1}\|=1,\ldots,\|\wtilde{Z}_{n,k}\|=1\}\big].
$$
Now suppose that $\E\big[(\det C_{n,k})^{-1/2} \big| \|\wtilde{Z}_{n,1}\|=1,\ldots,\|\wtilde{Z}_{n,k}\|=1\big]$ diverges.  Then we apply dilations to replace each $\wtilde{Z}_{n,j}$ by $s_j \wtilde{Z}_{n,j}$, where $s_1,\ldots,s_k > 0$.  Each vector $\wtilde{Z}_{n,j}$ remains normally distributed under these dilations, and $\det(C_{n,k})$ is transformed to $(s_1\cdots s_k)^2 \det(C_{n,k})$.  Therefore $\E\big[(\det C_{n,k})^{-1/2} \big| \{\|\wtilde{Z}_{n,1}\|=s_1,\ldots,\|\wtilde{Z}_{n,k}\|=s_k\}\big]$ diverges, for all $s_1,\ldots,s_k > 0$.  Integrating with respect to $s_1,\ldots,s_k$, it follows that $\E\big[(\det C_{n,k})^{-1/2} \big| \{\|\wtilde{Z}_{n,1}\| \le s_1,\ldots,\|\wtilde{Z}_{n,k}\| \le s_k\}\big]$ also diverges, and letting $s_1,\ldots,s_k \to \infty$ we deduce that the unconditional expectation, $\E[(\det C_{n,k})^{-1/2}]$, diverges.  

By Example \ref{ex_Xn_gaussian_3}, $\E[(\det C_{n,k})^{-1/2}] < \infty$ for $d_n \ge k+1$.  Therefore we deduce by contradiction that $\E[(\det B_{n,k})^{-1/2}] < \infty$ for all $n$ such that $d_n \ge k+1$, so \ref{condition_3} holds.  
}\end{example}

\begin{example}{\rm
\label{ex_spherical_Xn}
(Continuation of Example \ref{ex_Xn_ud_on_ball}): 
Suppose that $X_n$ is spherically distributed.  Then $X_n \eqdist R_n \Theta_n$ where $R_n \ge 0$, $\Theta_n$ is uniformly distributed on $\mcalS^{d_n-1}$, and $R_n$ and $\Theta_n$ are independent.  We assume that $\E(R_n^{-1}) < \infty$ for all $n$.  

Since $R_n$ and $\Theta_n$ are independent then 
$$
\varphi_{X_n}(t\theta) = \E_{R_n} \E_{\Theta_n} \exp(\i  R_n t\theta'\Theta),
$$
and by a change-of-variable, $t \to t/R_n$, we obtain 
\begin{align*}
\int_{-\infty}^\infty |\varphi_{X_n}(t\theta)| \dd t &= \int_{-\infty}^\infty |\E_{R_n} R_n^{-1} \E_{\Theta_n} \exp(\i  t\theta'\Theta)| \dd t 
= \E(R_n^{-1}) \! \int_{-\infty}^\infty |\E_{\Theta_n} \exp(\i  t\theta'\Theta)| \dd t.
\end{align*}
The latter integral is finite, as shown in Example \ref{ex_Bingham_contd_2}, and by assumption, $\E(R_n^{-1}) < \infty$, so it follows that \eqref{eq_integrable_cf_Yn_Xin} holds.  

Let $\wtilX_{n,1},\ldots,\wtilX_{n,k}$ be independent copies of $X_n$, with corresponding polar coordinates decompositions $\wtilX_{n,j} \eqdist \wtilde{R}_{n,j} \wtilde{\Theta}_{n,j}$, $j=1,\ldots,k$, and $\wtilde{R}_{n,1},\ldots,\wtilde{R}_{n,k},\wtilde{\Theta}_{n,1},\ldots,\wtilde{\Theta}_{n,k}$ are mutually independent.  Letting $B_{n,k} = \big(\wtilde{\Theta}_{n,j}'\wtilde{\Theta}_{n,r}\big)_{j,r=1}^k$, we obtain 
$$
\det (A_{n,k}) 
= \det\big(\wtilde{R}_{n,j}\wtilde{R}_{n,r} \wtilde{\Theta}_{n,j}'\wtilde{\Theta}_{n,r}\big)_{j,r=1}^k
= \Big(\prod_{j=1}^k \wtilde{R}_{n,j}^{~2}\Big) B_{n,k},
$$
and since $\wtilde{R}_{n,1},\ldots,\wtilde{R}_{n,k}$ are independent and identically distributed, and independent of $B_{n,k}$, then 
$$
\E[(\det A_{n,k})^{-1/2}] 
= \big(\E(R_n^{-1})\big)^k \E [(\det B_{n,k})^{-1/2}].
$$
By Example \ref{ex_Bingham_contd_2}, $\E [(\det B_{n,k})^{-1/2}] < \infty$ for $d_n \ge k+1$; also $\E(R_n^{-1}) < \infty$, by assumption.  Therefore $\E[(\det A_{n,k})^{-1/2}] < \infty$ for all $d_n \ge k+1$, so \ref{condition_3} holds.
}\end{example}

\begin{example}{\rm
\label{ex_Xn_ud_hyper_rect_contd_2}
(Continuation of Example \ref{ex_Xn_ud_hyper_rect}): 
We again assume that $L_n$, the length of each side of the hypercube $\mcalC^{d_n}(L_n)$, satisfies $d_n L_n^2 \cip 12\sigma^2$ as $n \to \infty$, and hence $(d_n L_n^2)^{-1/2} \cid (12\sigma^2)^{-1/2}$.  So we assume that $\E_{L_n} [(d_n L_n^2)^{-k/2}] < \infty$, and therefore $\E_{L_n} [(d_n L_n^2)^{-1/2}] < \infty$, for all sufficiently large $n$.  

Since $X_n = (X_{n;1},\ldots,X_{n;d_n})'$, conditional on $L_n$, is uniformly distributed on $\mcalC^{d_n}(L_n)$ then $X_{n;1}|L,\ldots,X_{n;d_n}|L$ are mutually independent and each uniformly distributed on the interval $[-L_n/2,L_n/2]$.  Using the well-known notation $\sinc t = (\sin t)/t$ if $t \neq 0$, and $\sinc 0 = 1$, we obtain, for $t \in \R$ and $\theta = (\theta_1,\ldots,\theta_{d_n})' \in \mcalS^{d_n-1}$, 
\begin{align}
\label{eq_sinc_product}
\varphi_{X_n|L_n}(t\theta) &= \E_{X_n|L_n} \exp(\i t\theta'X_n) \nonumber \\
&= \prod_{j=1}^{d_n} \E_{X_{n;j}|L_n} \exp(\i t\theta_j X_{n;j}) 
= \prod_{j=1}^{d_n} \sinc (\tfrac12 L_n \theta_j t).
\end{align}
Suppose that $\theta_1,\ldots,\theta_{d_n} \neq 0$, then by applying to \eqref{eq_sinc_product} the generalized H\"older inequality, we find that 
$$
\int_{-\infty}^\infty |\varphi_{X_n|L_n}(t\theta)| \dd t 
\le \bigg(\prod_{j=1}^{d_n} \int_{-\infty}^\infty |\sinc (\tfrac12 L_n |\theta_j| t)|^{d_n} \dd t\bigg)^{1/d_n}.
$$
Making the change-of-variable $t \to 2t/L_n |\theta_j|$ in the $j$th integral and simplifying the resulting product, we obtain 
\begin{equation}
\label{eq_cf_Xn_Ln_bound}
\int_{-\infty}^\infty |\varphi_{X_n|L_n}(t\theta)| \dd t 
\le 2 L_n^{-1} \bigg(\prod_{j=1}^{d_n} |\theta_j|\bigg)^{-1/d_n} \int_{-\infty}^\infty |\sinc t|^{d_n} \dd t.
\end{equation}
\cite{Borwein_etal}, during the proof of their Lemma 2, showed that there exists a universal constant $c_0$ such that  
$$
\int_{-\infty}^\infty |\sinc t|^{d_n} \dd t \le c_0 \, d_n^{-1/2}
$$
for all $d_n \ge 2$.  Therefore it follows from \eqref{eq_cf_Xn_Ln_bound} that 
\begin{equation}
\label{eq_cf_Xn_Ln_bound_2}
\E_{L_n} \int_{-\infty}^\infty |\varphi_{X_n|L_n}(t\theta)| \dd t \le 2 c_0 \, \E_{L_n}[(d_n L_n^2)^{-1/2}] \cdot \bigg(\prod_{j=1}^{d_n} |\theta_j|\bigg)^{-1/d_n} < \infty.
\end{equation}
Since 
$$
\E_{L_n} |\varphi_{X_n|L_n}(t\theta)| \ge |\E_{L_n} \varphi_{X_n|L_n}(t\theta)| = |\varphi_{X_n}(t\theta)|
$$
then, starting from the right-hand side \eqref{eq_cf_Xn_Ln_bound_2} and applying Fubini's theorem to interchange the integral and expectation, we obtain 
$$
\infty > \E_{L_n} \int_{-\infty}^\infty |\varphi_{X_n|L_n}(t\theta)| \dd t 
= \int_{-\infty}^\infty \E_{L_n} |\varphi_{X_n|L_n}(t\theta)| \dd t \ge \int_{-\infty}^\infty |\varphi_{X_n}(t\theta)| \dd t.
$$
Therefore \eqref{eq_integrable_cf_Yn_Xin} holds.  

Next let $\wtilX_{n,1},\ldots,\wtilX_{n,k}$ be mutually independent copies of $X_n$.  Conditional on $L_n$, the vectors $\wtilX_{n,1},\ldots,\wtilX_{n,k}$ are independent and uniformly distributed on $\mcalC^{d_n}(L_n)$.  Since $\mcalC^{d_n}(L_n) \subset \mcalB^{d_n}(R_n)$, where $R_n = d_n^{1/2} L_n/2$, then 
\begin{align*}
\E{}_{\wtilde{\mcalX}_n|L_n} [(\det A_{n,k})^{-1/2}] 
&= L_n^{-d_n k} \int_{\mcalC^{d_n}(L_n)} \cdots \int_{\mcalC^{d_n}(L_n)} \big(\det(\wtilde{x{\hskip1pt}}_{n,j}'\wtilde{x}_{n,r})\big)^{-1/2} \prod_{j=1}^{k} \dd \wtilde{x}_{n,j} \\
&\le L_n^{-d_n k} \int_{\mcalB^{d_n}(R_n)} \cdots \int_{\mcalB^{d_n}(R_n)} \big(\det(\wtilde{x{\hskip1pt}}_{n,j}'\wtilde{x}_{n,r})\big)^{-1/2} \prod_{j=1}^{k} \dd \wtilde{x}_{n,j}.
\end{align*}
Let $\wtilde{\Theta}_{n,1},\ldots,\wtilde{\Theta}_{n,k}$ be mutually independent and uniformly distributed on the unit ball $\mcalB^{d_n}(1)$, and define $B_{n,k} = (\wtilde{\Theta}_{n,j}'\wtilde{\Theta}_{n,r})_{j,r=1}^k$.  Substituting $\wtilde{x}_{n,j} = R_n \wtilde{\theta}_{n,j}$, $j=1,\ldots,k$, and simplifying the resulting expression, we obtain 
\begin{align*}
\E{}_{\wtilde{\mcalX}_n|L_n} [(\det A_{n,k})^{-1/2}] 
&\le L_n^{-d_n k} R_n^{-k + d_n k} [\textrm{Vol}(\mcalB^{d_n}(1))]^k \\
& \quad \times \int_{\mcalB^{d_n}(1)} \cdots \int_{\mcalB^{d_n}(1)} \big(\det(\wtilde{\theta{\hskip1pt}}_{n,j}'\wtilde{\theta}_{n,r})\big)^{-1/2} \prod_{j=1}^{k} \frac{\dd \wtilde{\theta}_{n,j}}{\textrm{Vol}(\mcalB^{d_n}(1))} \\
&= d_n^{d_n/2} 2^{k - d_n k} \, [\textrm{Vol}(\mcalB^{d_n}(1))]^k \, (d_n L_n^2)^{-k/2} \, \E [(\det B_{n,k})^{-1/2}].
\end{align*}
Applying the law of total expectation, we obtain 
\begin{align*}
\E[(\det A_{n,k})^{-1/2}] &= \E_{L_n} \E_{\wtilde{\mcalX}_n|L_n} [(\det A_{n,k})^{-1/2}] \\
&\le d_n^{d_n/2} 2^{k - d_n k} \, [\textrm{Vol}(\mcalB^{d_n}(1))]^k \, \E [(d_n L_n^2)^{-k/2}] \, \E [(\det B_{n,k})^{-1/2}].
\end{align*}
By Example \ref{ex_spherical_Xn}, $\E [(\det B_{n,k})^{-1/2}] < \infty$ for all $d_n \ge k+1$.  Also, 
$\E [(d_n L_n^2)^{-k/2}] < \infty$ for all sufficiently large $n$.  Therefore $\E[(\det A_{n,k})^{-1/2}] < \infty$ for all sufficiently large $n$, so \ref{condition_3} holds.  
}\end{example}

Finally we provide three examples of $\Xi_n$ for which \eqref{eq_Xi_n_cf} and \eqref{eq_G_dn_integrable}, the assumptions in Theorem \ref{thm_Lp_conv_cond_pdf_spher}, are valid.  In each example we have $\Xi_n = V Z_n$ where $V > 0$, $Z_n \sim \mcalN_{d_n}(0,I_{d_n})$, and $V$ and $Z_n$ are independent.  Therefore \eqref{eq_Xi_n_cf} holds for each example, so it remains only to verify \eqref{eq_G_dn_integrable}.

\begin{example}{\rm 
\label{eq_Xi_n_inv_chi_sq}
(i) Let $Q_\nu \sim \chi^2_\nu$ and $G$ be the distribution function of $V = (Q_\nu/\nu)^{-1/2}$.  
As noted in an earlier example, $\Xi_n$ has a multivariate $t$-distribution with $\nu$ degrees-of-freedom.  Also \eqref{eq_G_dn_integrable} holds since, for all $k=1,2,3,\ldots$,  
$$
\int_0^\infty v^{-k} \dd G(v) = \E(V^{-k}) 
= (\nu/2)^{-k/2} \, \frac{\Gamma((\nu+k)/2)}{\Gamma(\nu/2)}.
$$

(ii) For $\nu \ge 2$, let $V = (Q_\nu/\nu)^{1/2}$, so that $\Xi_n$ has a spherically symmetric multivariate Laplace distribution.  Then \eqref{eq_G_dn_integrable} holds since, for $k=1,\ldots,\nu-1$, 
$$
E(V^{-k}) = (\nu/2)^{k/2} \, \frac{\Gamma((\nu-k)/2)}{\Gamma(\nu/2)}.
$$
(iii) Let $V_0$ be a positive stable random variable with index $\alpha \in (0,1)$ and Laplace transform $\E \exp(-t V_0) = \exp(- 2^\alpha t^\alpha)$, $t \ge 0$.  Setting $V = V_0^{1/2}$, it is simple to show that $\Xi_n = V Z_n$ has a spherically symmetric stable distribution such that 
$\E \exp(\i u'\Xi_n) = \exp(- \|u\|^{2\alpha})$, $u \in \R^{d_n}$.  As shown by \citet{Brockwell_Brown}, 
$$
\E(V^{-k}) = \E(V_0^{-k/2}) = \frac{2^{-k/2} \, \Gamma(1 + (k/2\alpha))}{\Gamma(1+(k/2))},
$$
for all $k=1,2,3,\ldots$, and this moment formula also follows from a stochastic representation, established by \citet{Meintanis}, for $V_0$ in terms of the Weibull and exponential distributions.  Therefore \eqref{eq_G_dn_integrable} holds.  
}\end{example}

\section{Concluding remarks and open problems}
\label{sec_conclusions}

In this article, we have developed a general theory for random modulation of data $X_n$, under mild assumptions on the asymptotic properties of the data.  The theory is valid for the broad class of spherically symmetric distributions for the modulators $\Xi_n$, and we have illustrated the results by providing numerous examples.  The results derived here motivate several open questions, some of which we now describe.   

It was shown in Remark \ref{rem_non_Gaussian_Xin} that the proof of Theorem \ref{thm_resampling_Xi_n_X_n} does not hold for the case in which $\Xi_n$ has a spherically symmetric stable (non-Gaussian) distribution.  This raises the issue of whether there exists a variant of Theorem \ref{thm_resampling_Xi_n_X_n} for such $\Xi_n$.  

It would also be interesting to develop a theory of random modulation based on modulators whose distributions are neither Gaussian nor spherically symmetric.  In particular, a theory of random modulation in which at least one of $X_n$ and $\Xi_n$ is discrete is likely to find a wide range of applications \citep{Bishop_etal}.   

In Theorem \ref{thm_resampling_Xi_n_X_n}, the number $l$ of resamplings of $\Xi_n$ is fixed.  A reviewer noted that \cite{Meckes12}, in a study of the asymptotic distributions of random projections, derived the optimal growth rate of $l$; in our notation, Meckes proved that $l$ cannot grow faster than $\log d_n/(\log \log d_n)$.  As noted by the reviewer, Meckes' results raise the problem, in the spherically symmetric setting, of obtaining the optimal growth rate in Theorem \ref{thm_resampling_Xi_n_X_n} if $l$ is allowed to grow with $n$.  

As noted by the associate editor, Theorems \ref{thm_Lp_conv_cond_pdf_spher} and \ref{thm_Lp_conv_cond_cdf_spher} remain valid for polynomials (in the density and distribution functions, respectively) of degree at most $k$.  This raises the problem of extending those results to smooth functions.  In such a problem, $k$ must increase with $d_n$, and we pose as an open problem the question of how quickly $k$, as a function of $d_n$, can be allowed to grow.  

\medskip

\noindent
\textbf{Acknowledgments}.  We are very grateful to Lutz Duembgen and Shyamal Peddada for comments and suggestions that helped us to improve an initial version of this article.  We also thank the reviewer and the associate editor for providing us with insightful comments.  

\phantomsection


\addcontentsline{toc}{section}{References}


\begin{thebibliography}{3}

\bibitem[Bagyan(2015)]{Bagyan}
Bagyan, A.~(2015).  \textsl{Central Limit Theorems for Randomly Modulated Sequences of Random Vectors with Resampling and Applications to Statistics}.  \href{https://etda.libraries.psu.edu/catalog/26606}{Doctoral dissertation}, Pennsylvania State University.

\bibitem[Bagyan and Richards(2024)]{Bagyan_Richards_24}
Bagyan, A., and Richards, D.~(2024).  Complete asymptotic expansions and the high-dimensional Bingham distributions.  \textit{TEST}, \textbf{33}, 540--563.

\bibitem[Bagyan and Richards(2026)]{Bagyan_Richards_26}
Bagyan, A., and Richards, D.~(2026).  Some extensions of P\'olya's characterization of the normal distribution.  \textit{In preparation}.

\bibitem[Bickel, Gur, and Nadler(2018)]{Bickel_etal}
Bickel, P. J., Kur, G., and Nadler, B.~(2018).  Projection pursuit in high dimensions.  \textsl{Proc. Natl. Acad. Sci. USA}, \textbf{115}, 9151--9156.

\bibitem[Bingham(1974)]{Bingham}
Bingham, C.~(1974).  An antipodally symmetric distribution on the sphere.  \textit{Ann. Statist.}, 2, 1201--1225.

\bibitem[Bingham(1985)]{Bingham85}
Bingham, C.~(1985).  Review of \textsl{Statistics on Spheres}, by G. S. Watson. \textit{Ann. Statist.}, \textbf{13}, 838--844.

\bibitem[Bishop, \textit{et al.}(1975)]{Bishop_etal}
Bishop, Y. M., Fienberg, S. E., and Holland, P. W.~(1975). \textsl{Discrete Multivariate Analysis: Theory and Practice}. Springer, New York.

\bibitem[Blaabjerg, \textit{et al.}(1997)]{Blaabjerg_etal}
Blaabjerg, F., Pedersen, J. K., and Thoegersen, P.~(1997).  Improved modulation techniques for PWM-VSI drives.  \textit{IEEE Trans. Industrial Electron.}, \textbf{44}, 87--95.

\bibitem[Black(1953)]{Black}
Black, H. S.~(1953).  \textsl{Modulation Theory}. van Nostrand, Princeton, NJ.

\bibitem[Bogachev(1998)]{Bogachev}
Bogachev, V. I.~(1998).  \textsl{Gaussian Measures}.  American Mathematical Society, Providence, RI.

\bibitem[Borwein, \textit{et al.}(2010)]{Borwein_etal}
Borwein, D., Borwein, J. M., and Leonard, I. E.~(2010).  $L_p$ norms and the $\sinc$ function.  \textit{Amer. Math. Monthly}, \textbf{117}, 528--539.

\bibitem[Brockwell and Brown(1978)]{Brockwell_Brown}
Brockwell, P. J., and Brown, B. M.~(1978).  Expansions for the positive stable laws.  \textit{Z. Wahrscheinlichkeitstheorie verw. Gebiete}, \textbf{45}, 213--224.

\bibitem[Chow and Teicher(1988)]{Chow_Teicher}
Chow, Y. S. and Teicher, H.~(1988). \textsl{Probability Theory: Independence, Interchangeability, Martingales}, second edition. Springer, New York. 

\bibitem[Cochenour, \textit{et al.}(2011)]{Cochenour_etal}
Cochenour, B., Mullen, L., and Muth, J.~(2011).  Modulated pulse laser with pseudorandom coding capabilities for underwater ranging, detection, and imaging.  \textit{Appl. Opt.}, \textbf{50}, 6168--6178.

\bibitem[Davidov and Peddada(2013)]{Davidov_Peddada}
Davidov, O., and Peddada, S.~(2013).  The linear stochastic order and directed inference for
multivariate ordered distributions.  \textit{Ann. Statist.}, \textbf{41}, 1--40. 

\bibitem[Diaconis and Freedman(1984)]{Diaconis_Freedman}
Diaconis, P., and Freedman, D.~(1984).  Asymptotics of graphical projection pursuit. \textit{Ann. Statist.}, \textbf{12}, 793--815.

\bibitem[Duembgen and Del Conte-Zerial(2013)]{Duembgen_etal}
Duembgen, L., and Del Conte-Zerial, P.~(2013).  On low-dimensional projections of high-dimensional distributions.  In:~\textsl{From Probability to Statistics and Back:~High-Dimensional Models and Processes--A Festschrift in Honor of Jon A. Wellner}.  \textit{IMS Collections}, \textbf{9}, 91--104. Institute of Mathematical Statistics, Hayward, CA.

\bibitem[Eaton(1981)]{Eaton}
Eaton, M. L.~(1981). On the projections of isotropic distributions. \textit{Ann. Statist.}, \textbf{9}, 391--400.

\bibitem[Freedman and Lane(1980)]{Freedman}
Freedman, D., and Lane, D.~(1980).  The empirical distribution of Fourier coefficients.  \textit{Ann. Statist.}, \textbf{8}, 1244--1251.

\bibitem[Horn and Johnson(2013)]{Horn}
Horn, R. A., and Johnson, C. R.~(2013).  \textsl{Matrix Analysis}, second edition.  Cambridge University Press, New York.

\bibitem[Huber(1985)]{Huber} 
Huber, P.~(1985).  Projection pursuit.  \textit{Ann. Statist.}, \textbf{13}, 435--475.

\bibitem[Hwang and Lee(2020)]{Hwang}
Hwang, I.-P., and Lee, C.-H.~(2020).  Mutual interferences of a true-random LiDAR with other LiDAR signals.  \textit{IEEE Access}, \textbf{8}, 124123--124133.

\bibitem[Kawata(1972)]{Kawata}
Kawata, T.~(1972).  \textsl{Fourier Analysis in Probability Theory}. Academic Press, New York.

\bibitem[Kume and Walker(2009)]{Kume_Walker}
Kume, A., and Walker, S. G.~(2009).  On the Fisher-Bingham distribution.  \textit{Statist. Comput.}, \textbf{19}, 167--172.

\bibitem[Li and Yin(2007)]{Li_Yin}
Li, B., and Yin, X.~(2007).  On surrogate dimension reduction for measurement error regression: an invariance law.  \textit{Ann. Statist.}, \textbf{35}, 2143--2172.

\bibitem[Lieb(1973)]{Lieb}
Lieb, E. H.~(1973).  Convex trace functions and the Wigner-Yanase-Dyson conjecture.  \textit{Adv. Math.}, \textbf{11}, 267--288.

\bibitem[Lok and Lehnert(1998)]{Lok_Lehnert}
Lok, T. M., and Lehnert, J. S.~(1998).  An asymptotic analysis of DS/SSMA communication systems with general linear modulation and error control coding.  \textit{IEEE Trans. Inform. Theory}, \textbf{44}, 870--881.

\bibitem[Loperfido(2020)]{Loperfido}
Loperfido, N.~(2020).  Kurtosis-based projection pursuit for outlier detection in financial time series.  \textit{European J. Finance}, \textbf{26}, 142--164.

\bibitem[Magnus and Neudecker(2019)]{Magnus_Neudecker}
Magnus, J. R., and Neudecker, H.~(2019). \textsl{Matrix Differential Calculus with Applications in Statistics and Econometrics}, third edition. Wiley, Hoboken, NJ.

\bibitem[Malley(1983)]{Malley}
Malley, J. D.~(1983).  Statistical and algebraic independence.  \textit{Ann. Statist.}, \textbf{11}, 341--345.

\bibitem[Meckes(2009)]{Meckes}
Meckes, E. S.~(2009).  Quantitative asymptotics of graphical projection pursuit.  \textit{Electron. Commun. Probab.}, \textbf{14}, 176--185.

\bibitem[Meckes(2012)]{Meckes12}
Meckes, E. S.~(2012).  Projections of probability distributions: A measure-theoretic Dvoretzky theorem.  In: \textsl{Geometric Aspects of Functional Analysis}, Lecture Notes in Mathematics, vol. 2050, pp. 317--326.  Springer, New York.

\bibitem[Meintanis(2007)]{Meintanis}
Meintanis, S. G.~(2007).  Testing for exponentiality against Weibull and Gamma decreasing hazard rate
alternatives.  \textit{Kybernetika}, \textbf{43}, 307--314.

\bibitem[Muirhead(1982)]{Muirhead}
Muirhead, R. J.~(1982).  \textsl{Aspects of Multivariate Statistical Theory}.  Wiley, New York.

\bibitem[Olver and Wong(2025)]{Olver_Wong}
Olver, F. W. J., and Wong, R.~(2025).  \href{https://dlmf.nist.gov/2}{Chapter 2: Aymptotic Approximations}.  In: \textsl{NIST Digital Library of Mathematical Functions}, Release 1.2.4 of 2025-03-15 (F.~W.~J. Olver, \textit{et al.}, eds.).  

\bibitem[Papoulis(1983)]{Papoulis}
Papoulis, A.~(1983). Random modulation:~A review. \textit{IEEE Trans. Acoustics, Speech, Signal Proc.}, \textbf{31}, 96--105.

\bibitem[P\'olya(1923)]{Polya}
P\'olya, G.~(1923).  Herleitung des Gau{\ss}schen Fehlergesetzes aus einer Funktionsalsgleichung.  \textit{Math. Zeitschrift}, \textbf{18}, 96--108.

\bibitem[Reeves(2017)]{Reeves}
Reeves, G.~(2017).  Conditional central limit theorems for Gaussian projections. In: \textit{2017 IEEE International Symposium on Information Theory}, 3045--3049. IEEE.

\bibitem[Ressel(1976)]{Ressel}
Ressel, P.~(1976).  A short proof of Schoenberg's theorem.  \textit{Proc. Amer. Math. Soc.}, \textbf{57}, 66--68.

\bibitem[Ross(2010)]{Ross}
Ross, S.~(2010).  \textsl{A First Course in Probability}, eighth edition.  Prentice Hall, Upper Saddle River, NJ.

\bibitem[Roy, \textit{et al.}(2019)]{Roy_etal}
Roy, D., \textit{et al.}~(2019).  RFAL: Adversarial learning for RF transmitter identification and classification.  \textit{IEEE Trans. Cognitive Commun. Networking}, \textbf{6}, 783--801.

\bibitem[Schoenberg(1938)]{Schoenberg}
Schoenberg, I. J.~(1938).  Metric spaces and completely monotone functions. \textit{Ann. Math.}, \textbf{39}, 811--841.

\bibitem[She, \textit{et al.}(2011)]{She_etal}
She, C.-Y., \textit{et al.}~(2011).  Mesopause-region temperature and wind measurements with pseudorandom modulation continuous-wave (PMCW) lidar at 589 nm. \textit{ Appl. Opt.}, \textbf{50}, 2916--2926.

\bibitem[Steerneman and van Perlo-ten Kleij(2005)]{Steerneman}
Steerneman, A.~G.~M., and van Perlo-ten Kleij, F.~(2005).  Spherical distributions: Schoenberg (1938) revisited.  \textit{Expo. Math.}, \textbf{23}, 281--287.

\bibitem[Stern(2016)]{Stern}
Stern, A.~(2016).  Special aspects of the application of compressive sensing in optical imaging and sensing: Challenges, some solutions, and open questions. In: \textsl{Optical Compressive Imaging} (A. Stern, ed.).  CRC Press, Boca Raton.

\bibitem[Tang and Clement(2010)]{Tang_Clement}
Tang, S. C., and Clement, G. T.~(2010).  Standing-wave suppression for transcranial ultrasound by random modulation.  \textit{IEEE Trans. Biomed. Eng.}, \textbf{57}, 203--205.

\bibitem[van Trees(2002)]{vanTrees}
van Trees, H. L.~(2002).  \textsl{Detection, Estimation, and Modulation Theory, Part IV: Optimum Array Processing}.  Wiley, New York.

\bibitem[Wee and Tatikonda(2023)]{Wee_Tatikonda}
Wee, T. L. H., and Tatikonda, S.~(2025).  Random projections beyond zero overlap.  \textit{Electron. J. Probab.}, \textbf{30}, 1--26.

\bibitem[Yang, \textit{et al.}(2015)]{Yang_etal}
Yang, Z., Li, C., Yu, M., Chen, F., and Wu, T.~(2015).  Compact 405-nm random-modulation continuous wave lidar for standoff biological warfare detection.  \textit{J. Appl. Remote Sens.}, \textbf{9}, 096042.

\bibitem[Yanushkevichius and Yanushkevichiene(2007)]{YY2}
Yanushkevichius, R., and Yanushkevichiene, O.~(2007). Stability of a characterization by the identical distribution of linear forms.  \textit{Statistics}, \textbf{41}, 345--362.

\bibitem[Zolotarev(1986)]{Zolotarev}
Zolotarev, V.~M.~(1986). \textsl{One-Dimensional Stable Distributions}.  Translations of Mathematical Monographs.  American Mathematical Society, Providence, RI.

\end{thebibliography}
\end{document}